\documentclass[11pt]{amsart}
\usepackage{fullpage}
\usepackage{amssymb}
\usepackage{amstext}
\usepackage{amsmath}
\usepackage{amsthm}
\usepackage{stmaryrd}
\usepackage{verbatim}
\usepackage{hyperref}
\newtheorem{theorem}{Theorem}%[section]
\newtheorem{corollary}[theorem]{Corollary}
\newtheorem{proposition}[theorem]{Proposition}
\newtheorem{lemma}[theorem]{Lemma}
\theoremstyle{remark}
\newtheorem{remark}{Remark}

\newcommand{\spt}{\operatorname{spt}}
\newcommand{\SL}{\operatorname{SL}}
\begin{document}
\title{$\ell$-adic properties of modular forms of half-integral weight with eta-quotient multipliers}

\author{Wei Wang}
\address{Department of Mathematics, Shaoxing University, Shaoxing 312000, China}
\email{weiwang\_math@163.com}

\begin{abstract}
We establish general criteria for prime-power congruences and coefficient recurrences for half-integral weight modular forms with multiplier systems induced by eta-quotients. The results concern weakly holomorphic forms of weight at least $3/2$ and harmonic weak Maass forms of weight $3/2$ whose shadows are multiples of holomorphic eta-quotients of weight $1/2$. The criteria combine a common quadratic-symbol condition on all polar coefficients with a numerical condition on cusp orders and dimensions of integral-weight cusp-form spaces. No bound on the polar orders is required. Their proofs use formulas for the action of Hecke operators at arbitrary cusps and a vanishing argument based on cusp orders. As illustrations of these general criteria, we obtain congruences for sums of core sizes, the second moment of the color difference in nine-colored generalized Frobenius partitions, and smallest-parts functions.
\end{abstract}

%\today
\maketitle
\section{Introduction and main results}\label{section_introduction}
Half-integral weight modular forms with eta-quotient multipliers arise naturally in the theory of partitions and $q$-series. An important class consists of quotients of Eisenstein series by eta-quotients. Their Fourier coefficients are closely related to partition functions and weighted partition statistics, and the action of Hecke operators provides a method for studying their congruences modulo powers of primes. 

One application concerns sums of core sizes. For an integer $t\ge2$, let $a_t(n)$ be the sum of the sizes of the $t$-cores of all partitions of $n$, following Rolen \cite{rolen2017t}. Here the $t$-core of a partition is obtained by successively removing rim hooks of length $t$. These sums are related to Eisenstein series by
\begin{equation}\label{equation_prime_level_eisenstein_family}
F_t(\tau)=\frac{t\bigl(tE_2(t\tau)-E_2(\tau)\bigr)}{\eta(\tau)}=\sum_{n\ge0}\bigl(24a_t(n)+(t-1)(24n+t)p(n)\bigr)q^{n-1/24},
\end{equation}
where $p(n)$ denotes the partition function, $E_2(\tau)=1-24\sum_{n\ge1}\sigma_1(n)q^n$, $q=e^{2\pi i\tau}$, and $\sigma_1(n)=\sum_{d\mid n}d$.
Our results give congruences modulo arbitrary powers of primes for $D_t(n)=24a_t(n)+(t^2-1)p(n)$ at every level $t>1$ for which $S_2(\Gamma_0(t))=\{0\}$, and also at $t=27,32,36,49,50$, with the additional prime restrictions specified in Corollary~\ref{corollary_core_size_congruences}. For the prime $5$, classical partition congruences also yield congruences for $a_t(n)$ itself. For example, when $t=3$, we obtain the prime-power family
\begin{equation}\label{equation_intro_five_power}
a_3\left(5^{2m+1}n+\frac{71\cdot5^{2m}+1}{24}\right)\equiv0\pmod{5^m}\qquad(m\ge1,\ n\ge0).
\end{equation}
See Corollaries \ref{corollary_core_size_congruences} and
\ref{corollary_uncorrected_core_congruences} for the general statements.

Congruences for half-integral weight modular forms with eta multipliers have been studied from several points of view. Ahlgren, Andersen, and Dicks \cite{ahlgren2024shimura} constructed a Shimura lift for cusp forms whose multipliers are powers of the Dedekind eta multiplier twisted by Dirichlet characters, and used it to obtain quadratic congruences modulo prime powers. Wang and Yang \cite{zbMATH07595066} and Williams \cite{williams2025hecke} established prime-power congruences for Eisenstein quotients, including $E_{2r}(\tau)/\eta(\tau)$. Williams also treats powers of $\eta$ in the denominator and admissible cases at the prime $3$. These results have applications to moments of ranks and cranks and to higher-order smallest-parts functions. Related congruences modulo primes were studied in \cite[Corollary 6]{wang2026explicit}. Ahlgren, Bringmann, and Lovejoy \cite{zbMATH05941220} established congruences modulo arbitrary powers of primes for several smallest-parts functions, extending earlier results of Ono \cite{zbMATH06133929}. Related Hecke relations were developed by Ahlgren and Kim \cite{zbMATH06342625}. 
 
The main results provide general criteria for prime-power congruences for classes of weakly holomorphic modular forms and harmonic weak Maass forms with eta-quotient multipliers. These criteria are formulated in terms of modularity, integrality, expansions at the cusps, and dimension inequalities. We illustrate the scope of these criteria through applications to partition statistics, obtaining congruences for sums of core sizes, a second moment of nine-colored generalized Frobenius partitions, and smallest-parts functions. 

To state the main results, we introduce some common notation. For an odd integer $\kappa$, let $\lambda=(\kappa-1)/2$. For even integers $k\ge2$, let $d_k(N)$  denote the dimension of the space of cusp forms $S_k(\Gamma_0(N))$. Let $N$ be a positive integer and let $g(\tau)=\prod_{\delta\mid N}\eta(\delta\tau)^{r_\delta}$ be an eta-quotient, with $r_\delta\in\mathbb Z$. Let
\[
r=\sum_{\delta\mid N}r_\delta,\qquad
R=\sum_{\delta\mid N}\delta r_\delta,\qquad
\Delta=\prod_{\delta\mid N}\delta^{r_\delta}.
\]
Thus $g$ has weight $r/2$. The eta-quotient $g$ induces a multiplier system $\nu_g$, defined by
\[
g(\gamma\tau)=\nu_g(\gamma)(c\tau+d)^{r/2}g(\tau),\qquad\gamma=\begin{pmatrix}a&b\\c&d\end{pmatrix}\in\Gamma_0(N).
\]
We use $M_k(\Gamma_0(N),\nu)$, $M_k^!(\Gamma_0(N),\nu)$, and $H_k^!(\Gamma_0(N),\nu)$ for holomorphic modular forms, weakly holomorphic modular forms, and harmonic weak Maass forms, respectively, of weight $k$ and multiplier $\nu$. For a harmonic weak Maass form $f$, $\xi_k f$ is called the shadow of $f$. The definitions and cusp conventions are fixed in Section~\ref{subsection_cusp_expansions}.

For a cusp $\mathfrak a$ of width $h_{\mathfrak a}$, use the parameter $q_{\mathfrak a}=e^{2\pi i\tau/h_{\mathfrak a}}$, and let $\rho_{\mathfrak a}=\operatorname{ord}_{\mathfrak a}(g)$ denote the initial exponent of $g$ in this parameter. For a representative $1/c$, with $c>0$, write $\rho_c=\rho_{1/c}$. The width of the cusp and the order of $g$ there are
\begin{equation}\label{equation_explicit_eta_orders}
h_c=\frac{N}{\gcd(c^2,N)},\qquad \rho_c=\frac{h_c}{24}\sum_{\delta\mid N}
\frac{\gcd(c,\delta)^2}{\delta}r_\delta.
\end{equation}
For the eta-quotient $g$, we define
\begin{equation}\label{equation_signed_eta_cusp_excess}
\mathcal W_N(g)=\sum_{\mathfrak a}\lfloor\rho_{\mathfrak a}\rfloor
-\#\{\mathfrak a:\rho_{\mathfrak a}\in\mathbb Z_{>0}\},
\end{equation}
where the sum and the count run over inequivalent cusps. These orders are computed by \eqref{equation_explicit_eta_orders}; each divisor $c\mid N$ occurs with multiplicity $\varphi(\gcd(c,N/c))$, where $\varphi$ is Euler's totient function. If $g$ is holomorphic, then
\[
\mathcal W_N(g)=\sum_{\mathfrak a}\max\{0,\lceil\rho_{\mathfrak a}\rceil-1\}.
\]

Throughout the main results, the prime $\ell$ satisfies
\begin{equation*}\tag{P}\label{equation_prime_hypotheses}
\ell\ge3,\qquad \ell\nmid N,\qquad\ell=3\ \Longrightarrow\ 3\mid R.
\end{equation*}
Let
\begin{equation}\label{equation_main_arithmetic_notation}
d=d_\ell=\begin{cases}24,&\ell\ge5,\\8,&\ell=3,\end{cases}.
\end{equation}
We also set
\begin{equation}\label{equation_uniform_local_indices}
\epsilon_\ell(x)=\left(\frac{(-1)^{\lambda+1}(d/8)\Delta x}{\ell}\right).
\end{equation}
Here $dR/24,d\rho_c\in\mathbb Z$.
For a form $f$ of weight $\kappa/2$, write $f_c=f|_{\kappa/2}\sigma_c$, where the scaling matrices are fixed in Section~\ref{section_preliminaries}. If $f$ is a harmonic weak Maass form, let $f_c^+$ denote its holomorphic part at the cusp $1/c$. For a weakly holomorphic form, we set $f_c^+=f_c$. For forms with multiplier $\overline\nu_g$, write
\[
f_c^+(\tau)=\sum_{n\gg-\infty}a_c(n)q_c^{n/d}.
\]
We show in Section \ref{subsection_eta_multiplier_properties} that $a_c(n)=0$ unless $n\equiv-d\rho_c\pmod d$ and at infinity, the corresponding coefficients $a(n)$ vanish unless $n\equiv-dR/24\pmod d$. Coefficients at nonintegral indices are also taken to be zero. We write $\mathcal P_c^-(f)$ for the principal part of $f_c^+$, i.e. the sum of the negative-exponent terms. All congruences are interpreted in $\mathcal O/\ell^m\mathcal O$, where $\mathcal O$ is a fixed local coefficient ring at $\ell$ as specified in Section~\ref{subsection_cusp_expansions}. We call a coefficient $\ell$-integral if it belongs to $\mathcal O$.

Our first theorem concerns weakly holomorphic forms of half-integral weight. 
\begin{theorem}\label{theorem_weakly_holomorphic}
Let $g\in M_{r/2}^!(\Gamma_0(N),\nu_g)$ be an eta-quotient,
and let $\kappa\ge3$ be odd. Suppose that $k=(\kappa+r)/2$ is an
even integer with $k\ge2$ and that 
\begin{equation}\label{equation_dimension_cusp_criterion}
\mathcal W_N(g)\ge d_k(N)+d_2(N)-\delta_{k,2},
\end{equation}
where $\delta_{k,2}$ is the Kronecker delta. Use the notation in \eqref{equation_main_arithmetic_notation} and \eqref{equation_uniform_local_indices}, fix a prime $\ell$ satisfying \eqref{equation_prime_hypotheses}, and fix $\chi(\ell)\in\{1,-1\}$. Let $f\in M_{\kappa/2}^!(\Gamma_0(N),\overline\nu_g)$, and write
\[
f(\tau)=\sum_{n\gg-\infty}a(n)q^{n/d}.
\]
Assume that $a(n)\in\mathcal O$ for every $n$ and that the following conditions hold:
\begin{enumerate}
\renewcommand{\theenumi}{W\arabic{enumi}}
\renewcommand{\labelenumi}{\textup{(\theenumi)}}
\item\label{hypothesis_weak_principal} At every cusp, each nonzero polar coefficient satisfies
\[
n<0,\quad a_c(n)\ne0\quad\Longrightarrow\quad
\epsilon_\ell(-h_cn)=\chi(\ell).
\]
\item\label{hypothesis_weak_constant} If $\rho_c\in\mathbb Z_{\le0}$, then $a_c(0)=0$.
\end{enumerate}
Then the following congruences hold:
\begin{enumerate}
\renewcommand{\theenumi}{\roman{enumi}}
\renewcommand{\labelenumi}{\textup{(\theenumi)}}
\item\label{conclusion_weak_zero} For every $m\ge1$, if $n\in\mathbb Z$ satisfies $\chi(\ell)\epsilon_\ell(-n)=1$, then
\[
a(\ell^{2m}n)\equiv0\pmod{\ell^{(\kappa-2)m}}.
\]
\item\label{conclusion_weak_recurrence} For every $j\ge1$ and every $n\in\mathbb Z$, we have
\begin{equation}\label{equation_main_unified_recurrence}
a(\ell^{j+2}n)\equiv\ell^{\lambda-1}\chi(\ell)a(\ell^jn)\pmod{\ell^{(\kappa-2)\lceil j/2\rceil}}.
\end{equation}
\end{enumerate}
\end{theorem}

At a cusp with $\rho_c\notin\mathbb Z$, the translation law forces $a_c(0)=0$. For quotients $E/g$ with $E\in M_k(\Gamma_0(N))$, Remark~\ref{remark_weak_quotient_criterion} explains how to check (\ref{hypothesis_weak_principal})--(\ref{hypothesis_weak_constant}) from the numerator.

When $g$ is holomorphic, condition \eqref{equation_dimension_cusp_criterion} is automatic if $d_k(N)=0$. For $\kappa\ge5$, condition \eqref{equation_dimension_cusp_criterion} is equivalent to
\[
d_2(N)=0,\qquad S_{\kappa/2}(\Gamma_0(N),\overline\nu_g)=\{0\},
\qquad \rho_{\mathfrak a}\notin\mathbb Z_{>0}\ \text{at every cusp};
\]
see Proposition~\ref{proposition_cusp_dimension_equivalence}. 

Our second theorem uses the same notation with $\kappa=3$ and $r=1$, so $\lambda=1$ and $k=2$. Here $g$ also specifies the shadow, and its Hecke eigenvalue determines $\chi(\ell)$. 

\begin{theorem}\label{theorem_harmonic_weak_maass}
Let $g\in M_{1/2}(\Gamma_0(N),\nu_g)$ be a holomorphic eta-quotient such that
\begin{equation}\label{equation_main_cusp_inequality}
\mathcal W_N(g)\ge2d_2(N)-1.
\end{equation}
Use the notation in \eqref{equation_main_arithmetic_notation} and \eqref{equation_uniform_local_indices} with $\kappa=3$, fix a prime $\ell$ satisfying \eqref{equation_prime_hypotheses}, and put
\[
\chi(\ell)=
\begin{cases}
\epsilon_\ell(dR/24),&R\ne0,\\
1,&R=0.
\end{cases}
\]
Let $f\in H_{3/2}^!(\Gamma_0(N),\overline\nu_g)$, and write
\[
f^+(\tau)=\sum_{n\gg-\infty}a(n)q^{n/d}.
\]
Assume that $f_c^+$ has coefficients in $\mathcal O$ at every cusp
and that the following conditions hold:
\begin{enumerate}
\renewcommand{\theenumi}{H\arabic{enumi}}
\renewcommand{\labelenumi}{\textup{(\theenumi)}}
\item\label{hypothesis_harmonic_principal} At every cusp, each nonzero polar coefficient satisfies
\[
n<0,\quad a_c(n)\ne0\quad\Longrightarrow\quad
\epsilon_\ell(-h_cn)=\chi(\ell).
\]
\item\label{hypothesis_harmonic_shadow} The shadow satisfies $\xi_{3/2}(f)=c_0g$ for a constant
$c_0$, which may be zero.
\item\label{hypothesis_harmonic_stability} For every divisor $c\mid N$ with $\rho_c=0$, we have
\begin{equation}\label{equation_main_cusp_stability}
\ell^2\equiv1\pmod{\gcd(c,N/c)}.
\end{equation}
\end{enumerate}
Then the following congruences hold:
\begin{enumerate}
\renewcommand{\theenumi}{\roman{enumi}}
\renewcommand{\labelenumi}{\textup{(\theenumi)}}
\item\label{conclusion_harmonic_zero} For every $m\ge1$, if $n\in\mathbb Z$ satisfies $\chi(\ell)\epsilon_\ell(-n)=1$, then
\[
a(\ell^{2m}n)\equiv0\pmod{\ell^m}.
\]
\item\label{conclusion_harmonic_recurrence} For every $j\ge1$ and every $n\in\mathbb Z$, we have
\begin{equation}\label{equation_harmonic_unified_recurrence}
a(\ell^{j+2}n)\equiv\chi(\ell)a(\ell^jn)\pmod{\ell^{\lceil j/2\rceil}}.
\end{equation}
\end{enumerate}
\end{theorem}

If the only possible polar term at each cusp is $a_c(-d\rho_c)q_c^{-\rho_c}$, then (\ref{hypothesis_harmonic_principal}) follows from Lemma~\ref{lemma_eta_local_character}.

\begin{remark}\label{remark_harmonic_constant_alternative}
Hypothesis (\ref{hypothesis_harmonic_stability}) can be replaced by
\begin{equation}\label{equation_harmonic_constant_alternative}
a_{\ell^2c}(0)=\chi(\ell)a_c(0)\qquad(\rho_c=0),
\end{equation}
for every representative $1/c$, with the cusp expansions fixed above. In particular, (\ref{hypothesis_harmonic_stability}) is unnecessary if these constant terms all vanish. The original hypothesis is a sufficient condition depending only on $N$, $g$, and $\ell$, and allows nonzero constant terms. The proof is given after that of Theorem~\ref{theorem_harmonic_weak_maass}.
\end{remark}

The following two corollaries give convenient consequences when $d_2(N)=0$.
\begin{corollary}\label{corollary_weight_two_dimension_zero}
If $d_2(N)=0$, then \eqref{equation_main_cusp_inequality} and \eqref{equation_main_cusp_stability} hold automatically for every holomorphic eta-quotient $g$ of weight $1/2$ and every prime $\ell\ge3$ with $\ell\nmid N$, subject to $3\mid R$ when $\ell=3$. 
\end{corollary}

\begin{corollary}\label{corollary_eisenstein_over_eta}
Suppose that $d_2(N)=0$, and let
$g\in M_{1/2}(\Gamma_0(N),\nu_g)$ be a holomorphic eta-quotient such that
\[
0\leq\operatorname{ord}_{\mathfrak a}(g)\leq1
\]
at every cusp $\mathfrak a$. Let $E\in M_2(\Gamma_0(N))$, and fix a prime
$\ell\ge3$ with $\ell\nmid N$, and assume $3\mid R$ if $\ell=3$. Suppose that the Fourier coefficients of
$E$ at infinity are $\ell$-integral, and write
\[
\frac{E(\tau)}{g(\tau)}=\sum_{n\ge -dR/24}b(n)q^{n/d}.
\]
Then both congruences in Theorem~\ref{theorem_harmonic_weak_maass} hold with $a(n)$
replaced by $b(n)$ and with the same $R$, $\Delta$, and $\chi(\ell)$.
\end{corollary}

The coefficient congruences in our main theorems arise from a common
Hecke congruence. Under the respective hypotheses of
Theorems \ref{theorem_weakly_holomorphic} and \ref{theorem_harmonic_weak_maass}, we prove
\begin{equation}\label{equation_main_hecke_congruence}
T_{\ell^{2m}}f-\ell^{\lambda-1}\chi(\ell)T_{\ell^{2m-2}}f
\equiv0\pmod{\ell^{(\kappa-2)m}}\qquad(m\ge1)
\end{equation}
at every cusp. In the harmonic case, this congruence refers to the Fourier coefficients of the holomorphic parts of the displayed Hecke difference. Here $T_1=I$, $T_{\ell^2}=T(\ell^2)$, the higher operators are defined by \eqref{equation_hecke_recurrence}, and $\chi(\ell)\in\{1,-1\}$ is specified in each theorem. 

The formulas for the Fourier coefficients at the cusps in Theorem~\ref{theorem_hecke_operator} explain the quadratic-symbol conditions on the polar coefficients in the main results: these hypotheses make the polar coefficients of the Hecke difference \eqref{equation_main_hecke_congruence} divisible by the required prime power. Multiplication by $g$ then gives a weakly holomorphic form of even integral weight whose coefficients satisfy prescribed vanishing bounds at every cusp. Lemma~\ref{lemma_integral_cusp_vanishing} shows that the dimension inequalities in the main theorems force this product to vanish modulo that power. In the harmonic case, an additional correction first cancels the shadow, and the congruence in hypothesis (\ref{hypothesis_harmonic_stability}) controls constant terms at the cusps where $g$ has order zero. The coefficient congruences follow from the formula at infinity.

This paper is organized as follows.  Section~\ref{section_applications} presents the applications. Section~\ref{section_preliminaries} recalls the necessary background and records properties of eta-quotient multipliers, and Section~\ref{section_hecke_operators} derives the Hecke formulas at arbitrary cusps. Section~\ref{section_main_proofs} proves the vanishing criterion, the main theorems, and the alternative constant-term condition.  The proofs of applications and their calculations at the cusps are collected in Section~\ref{section_application_proofs}.

\section{Examples and applications}\label{section_applications}
Throughout this section, let $L_t(\tau)=tE_2(t\tau)-E_2(\tau)$. Partition statistics are taken to be zero at negative or nonintegral arguments, and $p(0)=1$. In the applications and their proofs, the congruences of the main theorems are expressed in the natural indices of $\sum_n b(n)q^{n-R/24}$ by identifying $b(n)=a\bigl(d(n-R/24)\bigr)$.

\subsection{Sums of core sizes}
We retain the notation $a_t(n)$, $F_t$, and $D_t(n)$ from the introduction. Identity \eqref{equation_prime_level_eisenstein_family} follows from the case $j=0$ of \cite[Theorem 1.1]{rolen2017t}. Related weighted partition sums were studied by Merca \cite{merca2024sums}. Rolen's Corollary 1.5 shows $a_t(n)\equiv np(n)\pmod{t^2}$;
we consider congruences modulo powers of primes not dividing $t$. For $t=2$, $a_2(n)=\operatorname{SOME}(n)$ is the difference between the sums of odd and even parts in all partitions of $n$. Bardhan and Saikia \cite[Theorem 3.1]{bardhan2026proofs} proved $a_2(n)\equiv0\pmod{5^\alpha}$ whenever $24n\equiv1\pmod{5^\alpha}$.

Let
\[
\mathcal T=\{2,3,4,5,6,7,8,9,10,12,13,16,18,25,27,32,36,49,50\}.
\]
The elements of $\mathcal T$ with $t\le25$ are precisely the levels $t>1$ with $d_2(t)=0$.

\begin{corollary}\label{corollary_core_size_congruences}\label{corollary_prime_level_eisenstein}
Let $t\in\mathcal T$, and let $\ell\ge5$ be a prime with $\ell\nmid t$.
For $t=36$, $49$, or $50$, impose the additional condition
\[
\begin{array}{c|ccc}
t&36&49&50\\\hline
\text{condition}&\left(\frac3\ell\right)=1&\ell\ne5&
\left(\frac{13}{\ell}\right)=1.
\end{array}
\]
Then both congruences in Theorem~\ref{theorem_harmonic_weak_maass} hold for
$D_t(n)$ with $R=\Delta=1$ and
$\chi(\ell)=\left(\frac{3}{\ell}\right)$. In particular, if $m\ge1$
and $\left(\frac{1-24n}{\ell}\right)=1$, then
\[
24a_t\left(\ell^{2m}n-\frac{\ell^{2m}-1}{24}\right)\equiv-(t^2-1)p\left(\ell^{2m}n-\frac{\ell^{2m}-1}{24}\right)\pmod{\ell^m}.
\]
\end{corollary}

For $\ell=5,7,11$, the partition term can be removed by using the classical congruences of Watson and Atkin.

\begin{corollary}\label{corollary_uncorrected_core_congruences}
Suppose that $t$ and $\ell$ satisfy the conditions of Corollary
\ref{corollary_core_size_congruences}, with $\ell\in\{5,7,11\}$.
Then both congruences in Theorem~\ref{theorem_harmonic_weak_maass} hold for
$a_t(n)$ with $R=\Delta=1$ and
$\chi(\ell)=\left(\frac3\ell\right)$. In particular, for $m\ge1$,
\[
a_t\left(\ell^{2m}n-\frac{\ell^{2m}-1}{24}\right)\equiv0\pmod{\ell^m}\quad\text{if }\left(\frac{1-24n}{\ell}\right)=1.
\]
\end{corollary}

\subsection{A second moment of Frobenius partitions}
This application uses a weight $6$ numerator with a nonzero component in $S_6(\Gamma_0(9))$, a space of dimension $3$; the numerator and the cusp expansions establishing this property are given in Section~\ref{subsection_proof_frobenius}.

Generalized Frobenius partitions were introduced by Andrews
\cite{andrews1984generalized}. For a nine-colored generalized Frobenius partition $\pi$, let $T(\pi)$ be the sum of the color labels in its top row minus the sum of those in
its bottom row, with the colors numbered $0,\ldots,8$.
Following Kolitsch \cite[Section 3]{kolitsch1992generalized}, let $c\phi_9(n,t)$
count the partitions of size $n$ for which $T(\pi)=t$.
We consider the second moment of this color difference,
\[
\mathcal M_{2,9}(n)=\sum_{t\in\mathbb Z}t^2c\phi_9(n,t).
\]
Exchanging the two rows changes the sign of $T$, so all odd moments vanish.

\begin{corollary}\label{corollary_frobenius_second_moment}
Let $\ell\ge5$ be a prime. For every $m\ge1$ and every integer $n$ with $\left(\frac{9-24n}{\ell}\right)=1$, we have
\[
\mathcal M_{2,9}\left(\ell^{2m}n-\frac{3(\ell^{2m}-1)}8\right)\equiv0\pmod{15\ell^m}.
\]
For every $j\ge1$ and every integer $u$ satisfying $\ell^ju+3\equiv0\pmod8$, we have
\[
\mathcal M_{2,9}\left(\frac{\ell^{j+2}u+3}{8}\right)
\equiv\left(\frac3\ell\right)\mathcal M_{2,9}\left(\frac{\ell^ju+3}{8}\right)
\pmod{15\ell^{\lceil j/2\rceil}}.
\]
\end{corollary}

\subsection{Further weakly holomorphic applications}
The first example has weight $7/2$ and gives congruences modulo $\ell^{5m}$. The second has weight $3/2$ and uses an eta denominator that is not holomorphic at every cusp.

\noindent\emph{Cubes of parts in overpartitions.}
Let $T_3(n)$ be the sum of the cubes of all parts in the overpartitions of $n$, with each part counted according to its multiplicity. Overlining does not affect the weight of a part, and we put $T_3(0)=0$. Write $\overline P(q)=\prod_{j\ge1}(1+q^j)/(1-q^j)=\eta(2\tau)/\eta(\tau)^2$.
Marking the occurrences of each part, we have
\begin{equation}\label{equation_overpartition_cubes_gf}
\begin{aligned}
\sum_{n\ge0}T_3(n)q^n
&=2\overline P(q)\sum_{j\ge1}\frac{j^3q^j}{1-q^{2j}}\\
&=\frac{(E_4(\tau)-E_4(2\tau))/120}{\eta(\tau)^2/\eta(2\tau)}.
\end{aligned}
\end{equation}
Here $E_4(\tau)=1+240\sum_{n\ge1}\sigma_3(n)q^n$, where $\sigma_3(n)=\sum_{d\mid n}d^3$.

\begin{corollary}\label{corollary_overpartition_cubes}
Let $\ell$ be an odd prime. For every $m\ge1$ and every integer $n$ with $\left(\frac{-n}{\ell}\right)=1$, we have
\[
T_3(\ell^{2m}n)\equiv0\pmod{\ell^{5m}}.
\]
For every $j\ge1$ and every integer $u$, we have
\[
T_3(\ell^{j+2}u)\equiv\ell^2T_3(\ell^ju)\pmod{\ell^{5\lceil j/2\rceil}}.
\]
\end{corollary}

\medskip\noindent\emph{Odd parts in colored partitions.}
We next consider partitions in which even parts have one color and odd parts have three colors, with unrestricted multiplicities. Let $c_3(n)$ count these partitions and let $O_3(n)$ be the sum of all their odd parts, counted with multiplicity. Thus
\begin{equation}\label{equation_odd_three_color_gf}
\begin{aligned}
P_3(q):=\sum_{n\ge0}c_3(n)q^n
&=\prod_{j\ge1}\frac1{(1-q^{2j})(1-q^{2j-1})^3},\\
\sum_{n\ge0}O_3(n)q^n
&=3P_3(q)\sum_{\substack{j\ge1\\j\text{ odd}}}\frac{jq^j}{1-q^j}.
 \end{aligned}
\end{equation}
The function $c_3(n)$ is the reciprocal crank parity function studied
by Amdeberhan and Merca \cite{amdeberhan2025crank} and Chen
\cite{chen2025congruences}. Define $b_3(n)=c_3(n)+8O_3(n)$. Its generating function is
\begin{equation}\label{equation_odd_three_color_quotient}
q^{1/24}\sum_{n\ge0}b_3(n)q^n=\frac{L_2(\tau)}{g_3(\tau)},\qquad g_3(\tau)=\frac{\eta(\tau)^3}{\eta(2\tau)^2}.
\end{equation}
Here $g_3$ has a pole at infinity.

\begin{corollary}\label{corollary_odd_three_color}
Let $\ell\ge5$ be prime and put $\chi=\left(\frac6\ell\right)$. For every $m\ge1$ and $n\ge0$ with $\left(\frac{-2(24n+1)}\ell\right)=1$, we have
\[
b_3\left(\ell^{2m}n+\frac{\ell^{2m}-1}{24}\right)\equiv0\pmod{\ell^m}.
\]
For every $j\ge1$ and every positive integer $u$ satisfying $\ell^ju\equiv1\pmod{24}$, we have
\[
b_3\left(\frac{\ell^{j+2}u-1}{24}\right)\equiv\chi b_3\left(\frac{\ell^ju-1}{24}\right)\pmod{\ell^{\lceil j/2\rceil}}.
\]
When $\ell=7$, both assertions also hold with $b_3$ replaced by $O_3$. In particular,
\[
O_3\left(7^{2m}n+\frac{7^{2m}-1}{24}\right)\equiv0\pmod{7^m}\qquad\text{if }n\equiv3,4,6\pmod7.
\]
\end{corollary}

\subsection{Applications of Theorem \ref{theorem_harmonic_weak_maass}}
Theorem~\ref{theorem_harmonic_weak_maass} also applies when the shadow is zero. The first example illustrates why its constant-term hypothesis differs from (\ref{hypothesis_weak_constant}): the generating function has a nonzero constant term at a cusp where its eta denominator has order zero. The two smallest-parts functions then illustrate nonzero shadows.

\noindent\emph{Odd parts without repetition.}
Let $r(n)$ count partitions of $n$ without repeated odd parts, and let $\operatorname{SODO}(n)$ be the sum of all odd parts in these partitions. Garvan and Sarma \cite[Theorem 1.12]{garvan2026sum} studied this statistic and proved congruences modulo $5$ and $7$. Define $d_O(n)=24\operatorname{SODO}(n)-r(n)$, where $r(0)=1$ and $\operatorname{SODO}(0)=0$, and set
\[
g_O(\tau)=\frac{\eta(\tau)\eta(4\tau)}{\eta(2\tau)},\qquad
Q(q)=\sum_{n\ge0}r(n)q^n=\frac{q^{1/8}}{g_O(\tau)}.
\]
The corrected coefficients have the generating function
\begin{equation}\label{equation_sodo_eisenstein}
\sum_{n\ge0}d_O(n)q^{n-1/8}=\frac{2L_2(\tau)-L_4(\tau)}{g_O(\tau)}.
\end{equation}

\begin{corollary}\label{corollary_sodo_congruences}
For every odd prime $\ell$, both congruences in Theorem \ref{theorem_harmonic_weak_maass} hold for $d_O(n)$ with $R=3$, $\Delta=2$, and $\chi(\ell)=\left(\frac2\ell\right)$. In particular, for $m\ge1$, if
$\left(\frac{1-8n}{\ell}\right)=1$ and $M=\ell^{2m}n-(\ell^{2m}-1)/8$, then
\[
24\operatorname{SODO}(M)\equiv r(M)\pmod{\ell^m}.
\]

\end{corollary}

\medskip\noindent\emph{Smallest-parts functions.}
We conclude with two applications with nonzero shadows. The function $\spt_\omega(n)$ counts smallest parts in partitions of $n$ in which every odd part is less than twice the smallest part. Its overpartition analogue $\overline{\spt}_\omega(n)$ has the additional condition that the smallest part is overlined;
see \cite{zbMATH06828531,zbMATH06533997}. Prime-power congruences for these functions were studied in \cite{zbMATH07567828,zbMATH08135742}. Their completions are obtained from the forms in \cite{ahlgren2015euler,zbMATH05941220}
by dilation and addition of Eisenstein quotients. Writing $(a;q)_j=\prod_{v=0}^{j-1}(1-aq^v)$, with the usual infinite-product convention, their generating functions are
\begin{align}
\sum_{n\ge1}\spt_\omega(n)q^n
&=\sum_{s\ge1}\frac{q^s}{(1-q^s)^2(q^{s+1};q)_s(q^{2s+2};q^2)_\infty},
\label{equation_sptomega_combinatorial_gf}\\
\sum_{n\ge1}\overline{\spt}_\omega(n)q^n
&=\sum_{s\ge1}\frac{q^s(-q^{s+1};q)_s(-q^{2s+2};q^2)_\infty}
{(1-q^s)^2(q^{s+1};q)_s(q^{2s+2};q^2)_\infty}.
\label{equation_barsptomega_combinatorial_gf}
\end{align}
These follow by fixing the smallest part $s$ and weighting its multiplicity. The corresponding completions and their expansions at the cusps are given in Section~\ref{subsection_proof_smallest_parts}. The first completion has a cuspidal shadow; the second has a shadow with a nonzero constant term at infinity, so it also uses hypothesis (\ref{hypothesis_harmonic_stability}).

\begin{corollary}\label{corollary_spt_omega}\label{corollary_bar_spt_omega}
For every prime $\ell\ge5$, both congruences in Theorem
\ref{theorem_harmonic_weak_maass} hold for the following coefficients and parameters:
\[
\begin{array}{c|c|c|c}
\text{sequence}&R&\Delta&\chi(\ell)\\\hline
\spt_\omega(n)&2&2&\left(\frac3\ell\right)\\
\overline{\spt}_\omega(n)-\frac18\overline p(n/2)&0&1&1.
\end{array}
\]
Here $\overline p(n/2)=0$ if $n$ is odd. All assertions for the second
row also hold for $\ell=3$. In particular, for $m\ge1$,
\[
\begin{aligned}
\spt_\omega\left(\ell^{2m}n-\frac{\ell^{2m}-1}{12}\right)
&\equiv0\pmod{\ell^m}
&&\text{if }\left(\frac{1-12n}{\ell}\right)=1,\\
\overline{\spt}_\omega(\ell^{2m}n)
&\equiv\frac18\overline p(\ell^{2m}n/2)\pmod{\ell^m}
&&\text{if }\left(\frac{-2n}{\ell}\right)=1.
\end{aligned}
\]
\end{corollary}

\section{Preliminaries and eta-quotient multipliers}\label{section_preliminaries}
\subsection{Modular forms and cusp expansions}\label{subsection_cusp_expansions}
We now fix the transformation and cusp conventions used in Theorems \ref{theorem_weakly_holomorphic} and \ref{theorem_harmonic_weak_maass}. For complex numbers $x, z \in \mathbb{C}$ with $z \neq 0$, we choose the principal branch of the logarithm to define $z^x = e^{x \log z}$, where $\log z = \log |z| + i \arg z$ with $-\pi < \arg z \le \pi$. For $\gamma=\begin{pmatrix}a&b\\c&d\end{pmatrix}\in\Gamma_0(N)$, we also write $\nu_g(a,b,c,d)=\nu_g(\gamma)$. If $\nu_\eta$ denotes the standard eta multiplier, from the matrix identity
\[
\begin{pmatrix} a & b\delta \\ c/\delta & d \end{pmatrix} = \begin{pmatrix} \delta & 0 \\ 0 & 1 \end{pmatrix} \begin{pmatrix} a & b \\ c & d \end{pmatrix} \begin{pmatrix} \delta & 0 \\ 0 & 1 \end{pmatrix}^{-1}
\]
we have the explicit formula
\[
\nu_g(a,b,c,d) = \prod_{\delta \mid N} \nu_\eta(a, b\delta, c/\delta, d)^{r_\delta}.
\]

Let $\kappa$ be an odd integer, and write $\lambda := \frac{\kappa-1}{2} \in \mathbb{Z}$. If $f(\tau)$ is a function on the upper half-plane $\mathbb{H}$ and $\gamma=\begin{pmatrix}
a &b\\
c &d 
\end{pmatrix}\in \mathrm{GL}^{+}_2(\mathbb{Q})$, we define
\[
(f|_{\kappa/2}\gamma)(\tau)=(ad-bc)^{\kappa/4}(c\tau+d)^{-\kappa/2}f\left(\frac{a\tau+b}{c\tau+d}\right).
\]
We denote by $H_{\kappa/2}^!(\Gamma_0(N), \nu_g)$ the space of harmonic weak Maass forms of weight $\kappa/2$ on $\Gamma_0(N)$ with the multiplier system $\nu_g$.
Namely, any $f(\tau) \in H_{\kappa/2}^!(\Gamma_0(N), \nu_g)$ is a smooth function on the upper half-plane $\mathbb{H}$ that is annihilated by the weight $\kappa/2$ hyperbolic Laplacian 
\[
\Delta_{\kappa/2} := -y^2 \left( \frac{\partial^2}{\partial x^2} + \frac{\partial^2}{\partial y^2} \right) + i \frac{\kappa}{2} y \left( \frac{\partial}{\partial x} + i \frac{\partial}{\partial y} \right),
\]
where $\tau=x+iy$, has at most linear exponential growth at the cusps, and satisfies the transformation law
\[
f\left(\frac{a\tau+b}{c\tau+d}\right) = \nu_g(a,b,c,d) (c\tau + d)^{\kappa/2} f(\tau),
\]
for all $\gamma = \begin{pmatrix} a & b \\ c & d \end{pmatrix} \in \Gamma_0(N)$. We refer the reader to \cite{zbMATH06828732} for the basic concepts of harmonic weak Maass forms. For a nonzero form to satisfy the transformation law at $-I\in\Gamma_0(N)$, we must have
\[
\kappa\equiv r\pmod 4.
\]

There is a complex anti-linear differential operator $\xi_{\kappa/2}$ defined by
\[
\xi_{\kappa/2}(f)(\tau) = 2i y^{\kappa/2} \overline{\frac{\partial}{\partial \bar{\tau}} f(\tau)}.
\]
Since
\[
\Delta_{\kappa/2}=-\xi_{2-\kappa/2}\circ\xi_{\kappa/2},
\] the image of a harmonic weak Maass form under $\xi_{\kappa/2}$ lies in $M_{2-\kappa/2}^!(\Gamma_0(N),\overline\nu_g)$. Its kernel is $M_{\kappa/2}^!(\Gamma_0(N),\nu_g)$, see \cite[Proposition 3.2]{bruinier2004two}. We call $\xi_{\kappa/2}(f)$ the shadow of $f$.

Every $f\in H_{\kappa/2}^!(\Gamma_0(N),\nu_g)$ has a Fourier expansion of the form below, with holomorphic part $f^+$ and non-holomorphic part $f^-$, see \cite[Section 3]{bruinier2004two}.
\begin{equation}\label{equation_fourier_expansion_f}
\begin{aligned} 
f(\tau) &= f^+(\tau) + f^-(\tau) \\ 
&= \sum_{n \gg -\infty} a^+(n) q^{n/24} + a^-(0)y^{1-\kappa/2} \\ 
&\quad + \sum_{n\ne0} a^-(n) \Gamma(1-\kappa/2, -4\pi ny/24) q^{n/24}, 
\end{aligned}
\end{equation}
where $\Gamma(a,x)$ is the incomplete gamma function. Only finitely many $a^-(n)$ with $n>0$ are nonzero. The translation law $f(\tau+1)=e^{2\pi iR/24}f(\tau)$ implies that $a^\pm(n)\ne0$ only if $n\equiv R\pmod{24}$.

It is well known that every cusp of $\Gamma_0(N)$ can be represented in the form $1/c$ for some positive integer $c$ satisfying $1\leq c\leq N$. We also allow arbitrary positive $c$ when comparing expansions at equivalent cusp representatives. More precisely, two cusps $a/c'$ and $1/c$, expressed as coprime fractions, are equivalent under $\Gamma_0(N)$ if and only if 
\begin{equation}\label{equation_cusp_condition}
\gcd(c', N) = \gcd(c, N) \quad \text{and} \quad c \equiv a c' \pmod{\gcd((c')^2, N)},
\end{equation}
(see \cite[\S 6.3]{cohen2017modular}).
For the cusp $1/c$, we fix the standard scaling matrix
\[
\sigma_{c} = \begin{pmatrix} 1 & 0 \\ c & 1 \end{pmatrix} \in \mathrm{SL}_2(\mathbb{Z}),
\]
and the associated cusp width $h_{c} = N / \gcd(c^2,N)$, and we define the parameter 
\[
R_c := \sum_{\delta\mid N}\frac{r_\delta\gcd(\delta,c)^2}{\delta}.
\]
For the same $f$, the Fourier expansion of $f(\tau)$ at the cusp $1/c$, with respect to the local uniformizer $q_{c} := e^{2\pi i \tau/h_c}$, similarly splits into $f_c(\tau) = f_c^+(\tau) + f_c^-(\tau)$, defined by
\begin{equation}\label{equation_fourier_cusp}
\begin{aligned}
f_c(\tau) &:= \left(f|_{\kappa/2}\sigma_c\right)(\tau)= f(\sigma_{c}\tau)(c\tau+1)^{-\kappa/2} \\
&= \sum_{n \gg -\infty} a_c^+(n) q_c^{n/24} + a_c^-(0)y^{1-\kappa/2} + \sum_{n\ne0} a_c^-(n) \Gamma(1-\kappa/2, -4\pi ny/(24h_c)) q_c^{n/24}.
\end{aligned}
\end{equation}
At each cusp only finitely many $a_c^-(n)$ with $n>0$ are nonzero. These coefficients vanish at every cusp if and only if the shadow is holomorphic at all cusps. Since the Fourier expansions of a modular form at equivalent cusps may differ by certain roots of unity, we adopt the convention throughout this paper that the Fourier coefficients of a modular form $f(\tau)$ at the cusp $1/c$ precisely refer to the Fourier coefficients of $f_c(\tau)$. After conjugation by $\sigma_c$, a generator of the stabilizer of $1/c$, modulo $\{\pm I\}$, acts by $\tau\mapsto\tau+h_c$. Under the translation $\tau \mapsto \tau + h_c$, we have 
\[
f_c(\tau+h_c) = e^{2\pi i h_c R_c/24} f_c(\tau).
\]
For a proof of this identity, see Lemma~\ref{lemma_trans_of_cuspc}.
Comparing coefficients, we obtain $a_c^\pm(n)\ne0$ only if $n\equiv h_cR_c\pmod{24}$. We denote by $\mathcal P_c^-(f)$ the principal part of the holomorphic part at $1/c$, thus $\mathcal P_c^-(f)=\sum_{n<0}a_c^+(n)q_c^{n/24}$. 

For the arithmetic statements, we fix a number field containing the holomorphic Fourier coefficients under consideration and the algebraic factors in the cusp transformations, together with a place above $\ell$. Let $\mathcal O$ be the ring of integers of the corresponding local field, enlarged by a finite extension if necessary. The term $\ell$-integral means belonging to $\mathcal O$, and all congruences modulo $\ell^m$ are interpreted in $\mathcal O/\ell^m\mathcal O$. For harmonic weak Maass forms, coefficient congruences concern the holomorphic parts of the expansions at the cusps.

\subsection{Properties of the eta multiplier}\label{subsection_eta_multiplier_properties}
For $\gamma = \begin{pmatrix} a & b \\ c & d \end{pmatrix} \in \mathrm{SL}_2(\mathbb{Z})$ with $c\ne0$, the multiplier system $\nu_{\eta}(a,b,c,d)=\nu_{\eta}(\gamma)$ is explicitly given by: 
\[
\nu_{\eta}(\gamma) = 
\begin{cases} 
\varepsilon(c,d) \left(\dfrac{d}{c}\right) e\left( \dfrac{(a+d)c - bd(c^2-1) - 3c}{24} \right) \quad &\text{if } c \text{ is odd,}\\
\left(\dfrac{c}{d}\right) e\left( \dfrac{(a+d)c - bd(c^2-1) + 3d - 3 - 3cd}{24} \right) \quad &\text{if } c \text{ is even,}
\end{cases}
\]
where $\left(\frac{\cdot}{\cdot}\right)$ denotes the standard Kronecker symbol, $e(x):=e^{2\pi ix}$, and the factor $\varepsilon(c,d)$ is defined by:
\[
\varepsilon(c,d) = 
\begin{cases} 
-1 & \text{if } c < 0 \text{ and } d < 0, \\ 
1 & \text{otherwise,}
\end{cases}
\]
For $c=0$, we have
\[
\nu_\eta(1,b,0,1)=e(b/24),\qquad
\nu_\eta(-1,b,0,-1)=-i e(-b/24).
\]
In particular, $\nu_\eta(-I)=-i$.
For a proof of this classical result, see \cite[Chapter 4, Theorem 2]{knopp2008modular}.
We define the exponent term $E(a,b,c,d)$ modulo 24 of the eta multiplier by
\[
E(a,b,c,d) = 
\begin{cases} 
(a+d)c - bd(c^2-1) - 3c \pmod{24}, \quad  &\text{if } c \text{ is odd,}\\
(a+d)c - bd(c^2-1) + 3d - 3 - 3cd \pmod{24},\quad &\text{if } c \text{ is even.}
 \end{cases}
\]

We require the following basic properties of the eta multiplier.
\begin{proposition}\label{proposition_easy_pro_eta}
Let $\gamma=\begin{pmatrix}a &b \\ c& d \end{pmatrix}\in\mathrm{SL}_2(\mathbb{Z})$. Then
\begin{enumerate}
\item  $\nu_{\eta}(-a, -b, -c, -d) = i\nu_{\eta}(a, b, c, d)$ if $c>0$.
\item  $\nu_{\eta}(d, -b, -c, a) = \nu_{\eta}(a, b, c, d)^{-1}$ if $c\neq 0$.
\item  $\nu_{\eta}(a, -b, -c, d) = \nu_{\eta}(a, b, c, d)^{-1}$ if $c\neq 0$.
\item  $\nu_{\eta}(d, b, c, a) = \nu_{\eta}(a, b, c, d)$ if $c\neq 0$.
\item $\nu_{\eta}(b,-a,d,-c)=\varepsilon_S(c,d)e(-1/8)\nu_{\eta}(a,b,c,d)$, where $\varepsilon_S(c,d)=-1$ if $c\geq 0$ and $d<0$ and $1$ otherwise.
\item $\nu_{\eta}(a+cm, b+dm, c, d) = e(m/24) \nu_{\eta}(a, b, c, d)$. 
\item $\nu_{\eta}(a,am+b,c,cm+d)=e(m/24)\nu_{\eta}(a, b, c, d)$.
\item  For $m\in\mathbb{Z}$, $\nu_{\eta}(1-mbd,mb^2,-md^2,1+mbd)= e(m/24)$.
\end{enumerate}
\end{proposition}
\begin{proof}
For (1), compare the automorphy factors using $(-c\tau-d)^{1/2}=-i(c\tau+d)^{1/2}$ when $c>0$. For (2), apply the transformation law to $\gamma^{-1}\gamma$ and use
\[
\left(\frac1{c\tau+d}\right)^{1/2}(c\tau+d)^{1/2}=1 \qquad(c\ne0).
\]
Property (3) follows by complex conjugation of the eta function, and (4) follows from (2) and (3).

Let $S=\left(\begin{smallmatrix}0&-1\\1&0\end{smallmatrix}\right)$ and $T=\left(\begin{smallmatrix}1&1\\0&1\end{smallmatrix}\right)$. For (5), compare the transformation laws for $\gamma S$ and $\gamma$, using $\nu_\eta(S)=e(-1/8)$. Evaluation of the principal square roots at $\tau=i$ gives the factor $-1$ precisely when $c\ge0$ and $d<0$. Properties (6) and (7) follow from $\eta(\tau+m)=e(m/24)\eta(\tau)$ by left and right translation, respectively.

For (8), let $M=\left(\begin{smallmatrix}1-mbd&mb^2\\-md^2&1+mbd\end{smallmatrix}\right)$. Since $\gcd(b,d)=1$, choose
$A=\left(\begin{smallmatrix}b&x\\d&y\end{smallmatrix}\right)\in\mathrm{SL}_2(\mathbb Z)$.
Then $M=AT^mA^{-1}$, and
\[
(-md^2(A\tau)+1+mbd)^{1/2}(d\tau+y)^{1/2}=(d(\tau+m)+y)^{1/2}.
\]
The transformation law shows $\nu_\eta(M)\nu_\eta(A)=\nu_\eta(AT^m)$, so (7) proves the assertion.
\end{proof}

\begin{proposition}\label{proposition_eta_multiplier_p_square}
For any prime $p\geq 5$, we have $\nu_{\eta}(a, b, p^2 c, d) = \nu_{\eta}(a, p^2 b, c, d)$ if $ad-p^2bc=1$ and $\nu_{\eta}(p^2a, b,  c, d) = \nu_{\eta}(a,  b, c, p^2d)$ if $p^2ad-bc=1$.
\end{proposition}
\begin{proof}
Suppose first that $c>0$. Since $\gcd(p,d)=1$, the Kronecker symbols in the two eta multipliers agree. By the explicit expression for $E(a,b,c,d)$, we obtain
\[
\begin{aligned}
&E(a,b,p^2c,d)-E(a,p^2b,c,d)\\
&\quad=(p^2-1)\bigl((a+d)c-bd(p^2c^2+1)-3c\bigr)\qquad(c\text{ odd}).
\end{aligned}
\]
For even $c$, the same expression holds with $-3c$ replaced by $-3cd$. Both differences are divisible by $24$, since $p^2\equiv1\pmod{24}$. The case $c=0$ follows directly from the explicit formulas for $\nu_\eta$, and $c<0$ follows from Proposition~\ref{proposition_easy_pro_eta} (1).

For the second assertion, apply Proposition~\ref{proposition_easy_pro_eta} (5) to both multipliers. Since $\varepsilon_S(c,p^2d)=\varepsilon_S(c,d)$, the result reduces to the first assertion.
\end{proof}

\begin{proposition}\label{proposition_etamul_squares}
Let $\gamma = \begin{pmatrix} a & b \\ c & d \end{pmatrix} \in \mathrm{SL}_2(\mathbb{Z})$. If $c \neq 0$, $a > 0$, and $a^2 \equiv 1 \pmod{24}$, then
\[
\nu_{\eta}(\gamma) = \left(\frac{2c}{a}\right) \varepsilon_a^{-1} e\left( \frac{a(b-c)}{24} \right).
\]
Here $\varepsilon_a$ is the fourth root of unity defined as $\varepsilon_a = 1$ if $a \equiv 1 \pmod 4$, and $\varepsilon_a = i$ if $a \equiv 3 \pmod 4$.
\end{proposition}

\begin{proof}
We first prove the following assertion: For $\begin{pmatrix} a & b \\ c & d \end{pmatrix} \in \mathrm{SL}_2(\mathbb{Z})$, if $c > 0$, $d > 0$, and $d^2 \equiv 1 \pmod{24}$, then
\begin{equation}\label{equation_eta_multiplier_d_square}
\nu_{\eta}(\gamma) = \left(\frac{2c}{d}\right) \varepsilon_d^{-1} e\left( \frac{d(b-c)}{24} \right).
\end{equation}
Since $ad - bc = 1$ and $d^2 \equiv 1 \pmod{24}$, we have
\begin{equation}\label{equation_amod24}
a \equiv d + bcd \pmod{24}.
\end{equation}
We evaluate $\nu_{\eta}(\gamma)$ based on the parity of $c$:

If $c$ is even, substituting \eqref{equation_amod24} into $E(a,b,c,d)$, we obtain
\begin{align*}
E(a,b,c,d)&= (2d + bcd)c - bd(c^2-1) + 3d - 3 - 3cd \\
&\equiv d(b-c) + 3d - 3 \pmod{24}.
\end{align*}
It can be verified directly that $e\left(\frac{3d-3}{24}\right) = \left(\frac{2}{d}\right) \varepsilon_d^{-1}$. Thus the eta multiplier simplifies to:
\[
\nu_{\eta}(\gamma) = \left(\frac{c}{d}\right) \left(\frac{2}{d}\right) \varepsilon_d^{-1} e\left( \frac{d(b-c)}{24} \right) = \left(\frac{2c}{d}\right) \varepsilon_d^{-1} e\left( \frac{d(b-c)}{24} \right).
\]

If $c$ is odd, the explicit formula gives $\left(\frac{d}{c}\right)$ and $E(a,b,c,d)$ is $(a+d)c - bd(c^2-1) - 3c$. Substituting \eqref{equation_amod24}, we find that $E(a,b,c,d)$ equals  $2cd + bd - 3c$ modulo 24. 
Applying quadratic reciprocity for positive odd integers $c$ and $d$, we have 
\[
\left(\frac{d}{c}\right) = \left(\frac{c}{d}\right) (-1)^{\frac{d-1}{2}\frac{c-1}{2}}= e((d-1)(c-1)/8)\left(\frac{c}{d}\right)
\]
and
\begin{align*}
3(d-1)(c-1)+E(a,b,c,d)&\equiv 5cd + bd - 3d - 6c + 3 \\
&\equiv bd - cd + 3d - 3 \pmod{24}.
\end{align*}
This completes the proof of \eqref{equation_eta_multiplier_d_square}.

If $c>0$ and $a>0$, apply \eqref{equation_eta_multiplier_d_square} to $\left(\begin{smallmatrix}d&b\\c&a\end{smallmatrix}\right)$ and use Proposition~\ref{proposition_easy_pro_eta} (4). This gives the required formula. If $c<0$, apply the positive case to $\left(\begin{smallmatrix}a&-b\\-c&d\end{smallmatrix}\right)$ and use Proposition~\ref{proposition_easy_pro_eta} (3). Since $(\frac{-1}{a})=\varepsilon_a^2$, we obtain
\begin{equation}\label{equation_asquare_cnegative}
\nu_\eta(\gamma)=\left(\frac{-2c}{a}\right)\varepsilon_a e\left(\frac{a(b-c)}{24}\right)=\left(\frac{2c}{a}\right)\varepsilon_a^{-1}e\left(\frac{a(b-c)}{24}\right),
\end{equation}
as asserted.
\end{proof}

\begin{lemma}\label{lemma_eta_cubes_three}
The two identities in Proposition~\ref{proposition_eta_multiplier_p_square} hold for $p=3$ after both sides are cubed. Moreover, if
$\gamma=\left(\begin{smallmatrix}a&b\\c&d\end{smallmatrix}\right)
\in\mathrm{SL}_2(\mathbb Z)$, $c\ne0$, and $a>0$ is odd, then
\[
\nu_\eta(\gamma)^3
=\left(\frac{2c}{a}\right)\varepsilon_a^{-3}e\left(\frac{a(b-c)}8\right).
\]
If $3\nmid N$ and $3\mid R$, then, for every $\gamma\in\Gamma_0(N)$
and every positive cusp denominator $c$,
\[
\nu_g(\gamma)^8=1,\qquad 3\mid h_cR_c,\qquad
\nu_g\left(\begin{smallmatrix}1&0\\8N&1\end{smallmatrix}\right)=1.
\]
\end{lemma}
\begin{proof}
In the proof of Proposition~\ref{proposition_eta_multiplier_p_square}, the Kronecker factors agree for every odd prime, and the exponent differences are divisible by $p^2-1$. For $p=3$ they are divisible by $8$, which proves the
cubed identities. For the second assertion, repeat the proof of
Proposition~\ref{proposition_etamul_squares} modulo $8$: an odd $d$ satisfies
$d^2\equiv1\pmod8$, so $a\equiv d+bcd\pmod8$. The even-$c$ computation implies $E\equiv d(b-c)+3d-3\pmod8$. For odd $c$, quadratic reciprocity adds $(d-1)(c-1)$ to the exponent modulo $8$, and $2(d-1)(c-1)\equiv0\pmod8$ implies the same result. Since $e((3d-3)/8)=(\frac2d)\varepsilon_d^{-3}$, interchange $a,d$ and use conjugation for $c<0$ exactly as above.

Put $R^*=N\sum_{\delta\mid N}r_\delta/\delta$. In $\mathbb Z_{(3)}$, every divisor $\delta$ of $N$ satisfies $\delta^{-1}\equiv\delta\pmod3$. Thus $R^*\equiv NR\equiv0\pmod3$ and $R_c\equiv R\equiv0\pmod3$; as $h_cR_c$ is an integer, it is divisible by $3$. For $c\ne0$, the explicit eta exponent, for either parity of $c/\delta$, satisfies
\[
E(a,b\delta,c/\delta,d)\equiv\frac{(a+d)c-bdc^2}{\delta}+bd\delta\pmod3.
\]
Its sum with coefficients $r_\delta$ is therefore zero modulo $3$. The signs and Kronecker factors disappear on taking eighth powers, proving $\nu_g(\gamma)^8=1$. The case $c=0$ follows from the translation formula. Finally,
$\nu_g(1,0,8N,1)=e(-R^*/3)=1$.
\end{proof}

We record the leading coefficients and integrality of the expansions of an eta-quotient at the cusps. The following lemma allows arbitrary integer eta exponents and does not require holomorphy at the cusps.

\begin{lemma}\label{lemma_cc_magnitude}
Let $g(\tau)=\prod_{\delta\mid N}\eta(\delta\tau)^{r_\delta}$ be an eta-quotient of level $N$, and let $C_c$ be its leading coefficient at $1/c$. Then
\[
C_c=\zeta_c\prod_{\delta\mid N}\left(\frac{\gcd(c,\delta)}{\delta}\right)^{r_\delta/2},\qquad
|C_c|^2=\prod_{\delta\mid N}\left(\frac{\gcd(c,\delta)}{\delta}\right)^{r_\delta},
\]
where $\zeta_c$ is a root of unity. Moreover,
\[
(g|_{r/2}\sigma_c)(\tau)=C_cq_c^{h_cR_c/24}U_c(q_c),\qquad U_c(t)\in1+t\mathbb Z[\zeta][[t]]
\]
for a root of unity $\zeta$ of order dividing $24N$. Thus $C_c$ is an $\ell$-adic unit and all Fourier coefficients of $g$ at $1/c$ are $\ell$-integral whenever $\ell\nmid N$.

If $\ell\nmid N$ and $f\in M_{\kappa/2}^!(\Gamma_0(N),\overline\nu_g)$ with $(\kappa+r)/2$ an even integer, then $\ell$-integrality of the Fourier coefficients of $f$ at infinity implies $\ell$-integrality at every cusp.
\end{lemma}
\begin{proof}
For each $\delta\mid N$, let $d_\delta=\gcd(c,\delta)$ and choose
$x,y\in\mathbb Z$ such that $(\delta/d_\delta)y-(c/d_\delta)x=1$.
Let
\[
\gamma_\delta=\begin{pmatrix}\delta/d_\delta&x\\c/d_\delta&y\end{pmatrix},
\qquad z'=\frac{d_\delta^2}{\delta}\tau-\frac{xd_\delta}{\delta}.
\]
Since $\gamma_\delta z'=\delta\tau/(c\tau+1)$, the  transformation of the eta function gives
\[
\eta\left(\frac{\delta\tau}{c\tau+1}\right)=\nu_\eta(\gamma_\delta)
\left(\frac{d_\delta}{\delta}\right)^{1/2}(c\tau+1)^{1/2}\eta(z').
\]
Since $\delta$ divides both $h_cc^2$ and $h_c\delta^2$, it divides $h_cd_\delta^2$. Put $m_\delta=h_cd_\delta^2/\delta\in\mathbb Z_{>0}$. The product expansion is
\[
\eta(z')=e\left(-\frac{xd_\delta}{24\delta}\right)q_c^{m_\delta/24}
\prod_{n\ge1}\left(1-e\left(-\frac{nxd_\delta}{\delta}\right)q_c^{nm_\delta}\right).
\]
Taking the product with exponents $r_\delta$, we obtain the stated expansion of $g$.

Now suppose that $f$ has $\ell$-integral coefficients at infinity. For a sufficiently large integer $A\ge0$, the form $G=\eta^{24A}gf$ is holomorphic of even integral weight on $\Gamma_0(N)$ with trivial multiplier. Its expansion at infinity is $\ell$-integral. Since $\ell\nmid N$, by the $q$-expansion principle, we see that the $q$-expansion of $G$ is $\ell$-integral at every cusp after adjoining the necessary roots of unity \cite[\S1.6]{katz1973p}. At $1/c$ the expansion of $\eta^{24A}g$ is $C_cq_c^{Ah_c+\rho_c}$ times a unit of $\mathcal O[[q_c]]$, and the proof is completed by dividing by $\eta^{24A}g$.
\end{proof}

\begin{lemma}\label{lemma_trans_of_cuspc}
Let $f$ have weight $\kappa/2$ and multiplier $\nu_g$ on $\Gamma_0(N)$, where $\kappa\equiv r\pmod4$. For every $c>0$,
\[
 f_c(\tau+h_c)=e(h_cR_c/24)f_c(\tau).
\]
\end{lemma}
\begin{proof}
The quotient $f/g$ has even integral weight $(\kappa-r)/2$ and trivial multiplier, so $f_c/g_c$ is $h_c$-periodic. By Lemma~\ref{lemma_cc_magnitude}, $g_c(\tau+h_c)=e(h_cR_c/24)g_c(\tau)$, which gives the formula.
\end{proof}

\section{Hecke operators at cusps}\label{section_hecke_operators}
Throughout this section $p\ge3$ is prime and $p\nmid N$; when $p=3$ we assume $3\mid R$. Put
\[
d_p=\begin{cases}24,&p\ge5,\\8,&p=3.\end{cases}
\]
For a rational number integral at $p$, its Legendre symbol is evaluated after reduction modulo $p$. The Hecke operator $T(p^2)$ is defined via the double coset
$\Gamma_0(N)\alpha\Gamma_0(N)$, where $\alpha=\left(\begin{smallmatrix}1&0\\0&p^2\end{smallmatrix}\right)$. The compatibility of $\nu_g$ with itself follows from \cite[Theorem~3.9]{zhou2023double}, applied with the same exponent vector on both sides and with double-coset parameter $p^2$. Indeed, the square condition in that theorem is automatic, and the remaining conditions are
\[
(p^2-1)R\equiv(p^2-1)R^*\equiv0\pmod{24},\qquad R^*=\sum_{\delta\mid N}\frac{N}{\delta}r_\delta.
\]
For $p\ge5$ these follow from $p^2\equiv1\pmod{24}$. For $p=3$, the condition $3\nmid N$ gives $R^*\equiv NR\pmod3$, so both congruences are equivalent to $3\mid R$. Thus our prime hypotheses are exactly the compatibility conditions. We use the double-coset construction of \cite[Lemma~3.1 and Proposition~3.2]{zhou2023double}, specifying its cocycle factors and normalization below before computing the action at arbitrary cusps.

The following theorem generalizes \cite[Proposition 4]{zbMATH00963239} to weakly holomorphic modular forms and harmonic weak Maass forms of half-integral weight with eta-quotient multipliers.
\begin{theorem}\label{theorem_hecke_operator}
Let $g(\tau)=\prod_{\delta\mid N}\eta(\delta\tau)^{r_{\delta}}$ and $\Delta = \prod_{\delta \mid N} \delta^{r_\delta}$. Assume \eqref{equation_prime_hypotheses} with $\ell=p$. Choose a positive integer $u_p$ such that $pu_p\equiv1\pmod{d_pN}$. There exists a Hecke operator $T(p^2)$ on the space of harmonic weak Maass forms $H_{\kappa/2}^!(\Gamma_0(N),\nu_g)$, where $\kappa = 2\lambda + 1$. For every $f\in H_{\kappa/2}^!(\Gamma_0(N),\nu_g)$, at each cusp, the holomorphic and non-holomorphic parts of the expansion of $T(p^2)f$ are obtained separately from the corresponding parts of the expansion of $f$.  Use the expansion \eqref{equation_fourier_cusp} with $24$ replaced by $d_p$, and denote its coefficients by $a_c^\pm(n)$. Then the corresponding coefficients $A_c^\pm(n)$ of $T(p^2)f$ satisfy
\begin{equation}\label{equation_local_hecke_coefficients}
A_c^\pm(n) = a_{p^2c}^\pm(p^2n)+\left(\frac{(-1)^\lambda(d_p/8)\Delta h_cn}{p}\right)p^{\lambda-1}a_c^\pm(n)+p^{2\lambda-1}a_{u_p^2c}^\pm(n/p^2)
\end{equation}
for all $n\in\mathbb Z$ in the holomorphic part and for $n\ne0$ in the non-holomorphic part, with coefficients at nonintegral indices interpreted as zero. The coefficient $A_c^-(0)$ satisfies
\begin{equation}\label{equation_local_hecke_nonholomorphic_constant}
A_c^-(0)=p^{\kappa-2}a_{p^2c}^-(0)+a_{u_p^2c}^-(0).
\end{equation}
The same holomorphic formula applies to weakly holomorphic forms. Moreover,
\begin{equation}\label{equation_hecke_xi_compatibility}
\xi_{\kappa/2}\bigl(T(p^2)f\bigr)=p^{\kappa-2}T(p^2)\bigl(\xi_{\kappa/2}f\bigr),
\end{equation}
where the operator on the right has weight $2-\kappa/2$ and multiplier $\overline\nu_g$. 
\end{theorem}

We first construct and normalize the operator and establish its compatibility with $\xi_{\kappa/2}$. We then prove the formulas for the Fourier coefficients at the cusps.

Write $\Gamma=\Gamma_0(N)$ and $\alpha=\left(\begin{smallmatrix}1&0\\0&p^2\end{smallmatrix}\right)$. For $A=\left(\begin{smallmatrix}a&b\\c&d\end{smallmatrix}\right) \in\mathrm{GL}_2^+(\mathbb Q)$, let
\[
j(A,\tau)=(\det A)^{-1/4}(c\tau+d)^{1/2},\qquad\omega(A,B)=\frac{j(A,B\tau)j(B,\tau)}{j(AB,\tau)}.
\]
With the principal square roots fixed above, $\omega(A,B)\in\{1,-1\}$ is independent of $\tau$ and satisfies the cocycle identity. In particular,
\[
(f|_{\kappa/2}A)|_{\kappa/2}B=\omega(A,B)^{-\kappa}f|_{\kappa/2}(AB),\qquad \nu_g(\gamma\gamma')=\omega(\gamma,\gamma')^\kappa\nu_g(\gamma)\nu_g(\gamma').
\]
For $\beta=\gamma_1\alpha\gamma_2\in\Gamma\alpha\Gamma$, define
\begin{equation}\label{equation_corrected_double_coset_multiplier}
\widetilde\nu_g(\beta)=\nu_g(\gamma_1)\nu_g(\gamma_2)\omega(\gamma_1,\alpha)^\kappa\omega(\gamma_1\alpha,\gamma_2)^\kappa.
\end{equation}
This is the extension in \cite[Lemma~3.1]{zhou2023double}. It is independent of
the factorization, and $\widetilde\nu_g(\alpha)=1$.

Define
\[
T(p^2)f=p^{\lambda-3/2}\sum_{\beta\in\Gamma\backslash\Gamma\alpha\Gamma}
\widetilde\nu_g(\beta)^{-1}f|_{\kappa/2}\beta.
\]
By \cite[Proposition~3.2]{zhou2023double}, the sum is finite, is independent of the coset representatives, and has the required transformation law. The transformation argument also applies to smooth functions. Since the weight Laplacian commutes with the slash action and the finite sum preserves the allowed cusp growth, $T(p^2)$ preserves $H_{\kappa/2}^!(\Gamma_0(N),\nu_g)$.
Finally, $\xi_{\kappa/2}(f|_{\kappa/2}\beta) =(\xi_{\kappa/2}f)|_{2-\kappa/2}\beta$, together with conjugation of the multiplier and the ratio $p^{\kappa/2-2}/p^{-\kappa/2}=p^{\kappa-2}$ of the normalizing factors, proves the stated identity with $\xi_{\kappa/2}$.

A set of representatives for $\Gamma_0(N)\setminus \Gamma_0(N)\alpha\Gamma_0(N)$ can be formed by the following $p^2 + p$ elements:
\begin{align*}
\alpha_b=\begin{pmatrix} 1 & b \\ 0 & p^2 \end{pmatrix}&=\begin{pmatrix} 1 & 0 \\ 0 & p^2 \end{pmatrix}\begin{pmatrix} 1 & b \\ 0 & 1 \end{pmatrix}~~(0\leq b<p^2),\\
\beta_h=\begin{pmatrix} p & h \\ 0 & p \end{pmatrix}&=\begin{pmatrix} 1 & 0 \\ p Ns & 1 \end{pmatrix}\begin{pmatrix} 1 & 0 \\ 0 & p^2 \end{pmatrix}\begin{pmatrix} p & h \\ -Ns & u \end{pmatrix}~~(0<h<p),\\
\sigma= \begin{pmatrix} p^2 & 0 \\ 0 & 1 \end{pmatrix}&= \begin{pmatrix} p^2 & -t \\ N & d \end{pmatrix}\begin{pmatrix} 1 & 0 \\ 0 & p^2 \end{pmatrix}\begin{pmatrix} p^2d & t \\ -N & 1 \end{pmatrix}.
\end{align*}
As for $\beta_h$ we choose and fix, for each $h$, two integers $s$ and $u$ so that $p u + Nsh = 1$. For $\sigma$, we choose two integers $t$ and $d$ so that $p^2d + Nt = 1$.
We use $\alpha_B$ and $\beta_H$ for the same matrices with arbitrary
integer upper-right entries. With our square-root convention, the cocycle factors in these factorizations are all $1$: the upper triangular factors have positive
diagonal entries, and the remaining square-root factors have opposite principal arguments.

Thus
\begin{align*}
(T(p^2)f)(\tau) &= p^{\lambda-\frac32} \sum_{0\leq b<p^2}\nu_g(1,b,0,1)^{-1}f((\tau+b)/p^2)p^{-\kappa/2}\\
&\quad+p^{\lambda-\frac32}\sum_{0<h<p}\nu_g(1,0,pNs,1)^{-1}\nu_g(p,h,-Ns,u)^{-1}f(\tau+h/p)\\
&\quad+p^{\lambda-\frac32}\nu_g(p^2,-t,N,d)^{-1}\nu_g(p^2d,t,-N,1)^{-1}f(p^2\tau)p^{\kappa/2}. 
\end{align*}

It follows immediately from the definition of $\nu_{g}$ that
\begin{equation}\label{equation_nu_g_first}
\nu_g(1,b,0,1)=\prod_{\delta\mid N}\nu_{\eta}(1,b\delta,0,1)^{r_{\delta}}=e\left(\frac{bR}{24}\right)
\end{equation}
and similarly 
\begin{equation}\label{equation_nu_g_pns}
\nu_g(1,0,pNs,1)=e\left(\frac{-pNs\sum_{\delta\mid N}r_{\delta}/\delta  }{24}\right).
\end{equation}
For $p\ge5$, use Proposition~\ref{proposition_etamul_squares}. For $p=3$, use its cubed version in Lemma~\ref{lemma_eta_cubes_three}, multiply over $\delta$, and recover the product identity from its cube. This is valid because both sides are eighth roots of unity: $3\mid R$ and $3\mid R^*$, where $R^*=N\sum r_\delta/\delta$. We obtain
\begin{align*}
\nu_g(p,h,-Ns,u)&=\prod_{\delta\mid N}\nu_{\eta}(p,h\delta,-Ns/\delta,u)^{r_{\delta}}\\
&=e\left(\frac{pNs\sum_{\delta\mid N}r_{\delta}/\delta}{24}\right)\left(\frac{(-1)^r 2Ns\Delta}{p}\right)\varepsilon_p^{-r}e\left(\frac{phR}{24}\right).
\end{align*}
Since $\kappa\equiv r\pmod{4}$ and $Nsh\equiv 1\pmod{p}$, combining with \eqref{equation_nu_g_pns}, we have
\begin{equation}\label{equation_nu_g_second}
\nu_g(1,0,pNs,1)\nu_g(p,h,-Ns,u)=\left(\frac{2\Delta h}{p}\right)\varepsilon_p^{\kappa}e\left(\frac{phR}{24}\right).
\end{equation}
Finally, by Proposition~\ref{proposition_eta_multiplier_p_square}, its cubed version when $p=3$, and Proposition~\ref{proposition_easy_pro_eta} (3),
\begin{equation}\label{equation_nu_g_third}
\nu_{g}(p^2,-t,N,d)=\nu_{g}(1,-t,N,p^2 d)=\nu_{g}(p^2d,t,-N,1)^{-1}.
\end{equation}

Combining \eqref{equation_nu_g_first}, \eqref{equation_nu_g_second} and \eqref{equation_nu_g_third}, we get the explicit expression for the Hecke operator:
\begin{align*}
(T(p^2)f)(\tau) &= p^{-2} \sum_{0\leq b<p^2}e\left(-\frac{bR}{24}\right)f((\tau+b)/p^2)\\
&\quad+p^{\lambda-3/2}\sum_{0<h<p}\left(\frac{2\Delta h}{p}\right)\varepsilon_p^{-\kappa}e\left(-\frac{phR}{24}\right) f(\tau+h/p)\\
&\quad+p^{2\lambda-1}f(p^2\tau).
\end{align*}

The construction above proves the existence of the operator and its compatibility with $\xi_{\kappa/2}$. To prove the remaining formulas for the Fourier coefficients at the cusps, we use the following three lemmas, whose proofs are given at the end of this section.
\begin{lemma}\label{lemma_etamul_gamma_b}
Let $b, c, B ,\delta$ be integers such that  $c, B, \delta> 0$ and $\delta\mid B$.  Define the  matrix $\gamma_{b,\delta}$ as
\[
\gamma_{b,\delta} = \begin{pmatrix} 1+(1+bc)Bc & \delta\left(b - B - bBc\right) \\ Bc^2/\delta & 1-Bc \end{pmatrix}.
\]
Then $\gamma_{b,\delta}\in \mathrm{SL}_2(\mathbb{Z})$, and 
\[
\nu_{\eta}(\gamma_{b,\delta}) = e\left( \frac{b\delta}{24} - \frac{B\gcd(c,\delta)^2}{24\delta} \right).
\]
\end{lemma}

\begin{lemma}\label{lemma_etamul_gamma_h}
Let $p\ge5$ be prime, and let $B,c,h,\delta$ be positive integers
such that
\[
24\delta\mid B,\qquad Bc\equiv1\pmod p,\qquad h(1-Bc)\equiv pB\pmod{p^2}.
\]
Then
\[
\gamma_{h,\delta}=
\begin{pmatrix}
p+pBc+Bc^2h & \dfrac{h(1-Bc)-pB}{p^2}\delta\\
pBc^2/\delta & \dfrac{1-Bc}{p}
\end{pmatrix}
\in\mathrm{SL}_2(\mathbb Z),
\]
and
\[
\nu_\eta(\gamma_{h,\delta})=\left(\frac{2h\delta}{p}\right)\varepsilon_p
e\left(\frac{hp\delta}{24}\right).
\]
For $p=3$, with $8\delta\mid B$ in place of $24\delta\mid B$, the same multiplier identity holds after both sides are cubed.
\end{lemma}

\begin{lemma}\label{lemma_etamul_xi_b}
Let $p\ge5$ be prime, and let $b,c,H,\delta$ be positive integers
such that
\[
p\mid 1+bc,\qquad p^2\nmid 1+bc,\qquad (1+bc)H\equiv bp\pmod{p^2},\qquad 24\delta\mid H.
\]
Then
\[
\xi_{b,\delta}=
\begin{pmatrix}
\dfrac{p+(1+bc)cH}{p^2}
&\dfrac{bp-(1+bc)H}{p^2}\delta\\
c^2H/\delta&p-cH
\end{pmatrix}
\in\mathrm{SL}_2(\mathbb Z),
\]
and
\[
\nu_\eta(\xi_{b,\delta})=\left(\frac{2\delta H}{p}\right)\varepsilon_p^{-1}
e\left(\frac{b\delta}{24}-\frac{pH\gcd(c,\delta)^2}{24\delta}\right).
\]
For $p=3$, with $8\delta\mid H$ in place of $24\delta\mid H$, the same multiplier identity holds after both sides are cubed.
\end{lemma}

\begin{proof}[Completion of the proof of Theorem~\ref{theorem_hecke_operator}]
Recall that $d_p\rho_c=d_ph_cR_c/24\in\mathbb Z$. For $p=3$, every identity for $\nu_g$ below follows by applying the cubed identities for $\nu_\eta$ and then taking products.
Both sides are eighth roots of unity by Lemma~\ref{lemma_eta_cubes_three} and $3\mid h_cR_c$. Thus equality of cubes implies the displayed identity; no identity for an individual eta multiplier is being assumed at $p=3$. We first reduce to the case $p\nmid c$. If $p\mid c$, put $\widetilde c=c+d_pN$. Then $p\nmid\widetilde c$ and
\[
\sigma_{\widetilde c}=u_{d_pN}\sigma_c,\qquad u_L=\begin{pmatrix}1&0\\L&1\end{pmatrix}.
\]
For every positive multiple $L$ of $d_pN$,
\[
\nu_g(u_L)=\prod_{\delta\mid N}e(-L/(24\delta))^{r_\delta}=1,
\]
where for $p=3$ we use $3\mid R^*$. The principal square roots introduce no sign in this product of lower triangular matrices with positive lower-left entries. Hence $f_{c+L}=f_c$. The same holds for $T(p^2)f$, and also for the pairs
$p^2\widetilde c,p^2c$ and $u_p^2\widetilde c,u_p^2c$. Their widths and $R_c$ are unchanged. Thus the formula for $\widetilde c$ gives the formula for $c$, and we may assume $p\nmid c$.

By linearity, we may compute the action of $T(p^2)$ on the holomorphic and non-holomorphic parts separately. Apart from the term $a_c^-(0)y^{1-\kappa/2}$, the non-holomorphic part at $1/c$ is a linear combination, for $n\ne0$, of
\[
\Gamma\left(1-\frac{\kappa}{2},-\frac{4\pi ny}{d_ph_c}\right)q_c^{n/d_p}.
\]
Under the transformations $\tau\mapsto (\tau+B)/p^2$ and $\tau\mapsto p^2\tau$, the imaginary part changes to $y/p^2$ and $p^2y$, respectively. These positive rescalings preserve the chosen branch, so the incomplete gamma factors transform in accordance with the corresponding powers of $q_c$ for either sign of $n$, and the same exponential sums and quadratic Gauss sums arise as in the holomorphic case. Thus the computation below applies to all holomorphic coefficients and to the non-holomorphic coefficients with $n\ne0$, and in these ranges we can write $a_c^\pm(n)$. The $n=0$ term in the non-holomorphic part is treated separately at the end.

The double-coset sum defining $T(p^2)$ has $p^2+p$ terms. We first collect the $p^2$ terms involving the expansion at the cusp $1/(p^2c)$, using the following partition of the residue classes $b'\pmod{p^2}$.

For the $p^2-p$ values of $b \in \{0, \dots,p^2-1\}$ satisfying $\gcd(1+bc, p) = 1$, we define $b' \in \{0, \dots, p^2-1\}$ as the unique solution to the congruence
\[
(1+bc)b' \equiv -b \pmod{p^2}.
\]
For any $b' \in \mathbb{Z}$, we define $B = b'(p^2u_p^2-1)$  and the matrix
\[
\gamma_b = \begin{pmatrix} 1 + (1+bc)Bc & \frac{b}{p^2} - \frac{(1+bc)B}{p^2 } \\ p^2 Bc^2 & 1 - Bc \end{pmatrix}.
\]
The congruence for $b'$ ensures the top-right entry is integral.
Since $pu_p \equiv 1 \pmod {N}$, we have $\gamma_b \in \Gamma_0(N)$.
Direct matrix multiplication verifies 
\[
\alpha_b\sigma_c=\gamma_b\sigma_{p^2c}\alpha_B,
\]
where the cocycle factors are $1$. For $B>0$, this follows by letting $\operatorname{Im}\tau\to\infty$: in the decomposition the factors from $\gamma_b$ and $\sigma_{p^2c}$ have limiting arguments $0$ and $\pi/2$, respectively. The case $b=B=0$ is the identity case. The same argument applies to $\gamma_h$, $\gamma_0$, and $\xi_b$ below. Thus
\begin{equation}\label{equation_hecke_gamma_b}
(f|_{\kappa/2}\alpha_b)|_{\kappa/2}\sigma_c =(f|_{\kappa/2} \gamma_b \sigma_{p^2c} \alpha_B)=\nu_g(\gamma_b)(f|_{\kappa/2}\sigma_{p^2c})|_{\kappa/2}\alpha_B.
\end{equation}
For $B>0$, apply Proposition~\ref{proposition_eta_multiplier_p_square} to each eta multiplier, moving the factor $p^2$ from the lower-left entry to the upper-right entry, and then apply Lemma~\ref{lemma_etamul_gamma_b}. The case $b=B=0$ is  immediately. In both cases, we obtain
\[
\nu_g(\gamma_b)=e\left(\frac{bR-BR_c}{24}\right).
\]

The remaining $p$ residue classes omitted by $\alpha_b$ satisfy $c b' \equiv -1 \pmod p$. If $c b'\not\equiv -1\pmod {p^2}$, for each $h \in \{1, \dots, p-1\}$, we define $b'$ via the congruence
\[
h(1 + c b') \equiv -pb' \pmod{p^2}.
\]
Write $B=b'(p^2u_p^2-1)$ and consider the corresponding matrix
\[
\gamma_h = \begin{pmatrix} p + (p+hc)Bc &  \frac{h-(p+hc)B}{p^2 } \\ p B c^2 & \frac{1-Bc}{p} \end{pmatrix}.
\]
The conditions $c b' \equiv -1 \pmod p$ and $pu_p\equiv 1\pmod N$ ensure  $\gamma_h \in \Gamma_0(N)$. We also have
\[
\beta_h\sigma_c=\gamma_h\sigma_{p^2c}\alpha_B,
\]
which means
\begin{equation}\label{equation_hecke_gamma_h}
(f|_{\kappa/2}\beta_h)|_{\kappa/2}\sigma_c=\nu_g(\gamma_h)(f\big|_{\kappa/2}\sigma_{p^2c})\big|_{\kappa/2}\alpha_B.
\end{equation}
Since $B=b'(p^2u_p^2-1)$ and $pu_p\equiv1\pmod{d_pN}$, we have $d_p\delta\mid B$ for every $\delta\mid N$. Lemma~\ref{lemma_etamul_gamma_h} therefore applies and gives:
\[
\nu_g(\gamma_h)=\left(\frac{2\Delta h}{p}\right)\varepsilon_p^{\kappa}e\left(\frac{phR-BR_c}{24}\right).
\]

The remaining residue class corresponds to $c b' \equiv -1 \pmod{p^2}$. We define
\[
\gamma_{0} = \begin{pmatrix} p^2(1+Bc) & -B \\ B c^2 & \frac{1-Bc}{p^2} \end{pmatrix},
\]
which is integral under $c b' \equiv -1 \pmod{p^2}$, satisfying 
\[
\sigma\sigma_c=\gamma_0\sigma_{p^2c}\alpha_B,
\]
and thus
\begin{equation}\label{equation_hecke_gamma_0}
(f|_{\kappa/2}\sigma)|_{\kappa/2}\sigma_c=\nu_g(\gamma_0)(f\big|_{\kappa/2}\sigma_{p^2c})\big|_{\kappa/2}\alpha_B.
\end{equation}
Using Proposition~\ref{proposition_eta_multiplier_p_square} and Lemma~\ref{lemma_etamul_gamma_b}, we have
\[
\nu_g(\gamma_0)=\nu_g(1+Bc,-B,Bc^2,1-Bc)=e\left(\frac{-BR_c}{24}\right).
\]

The unified action over the $p^2$ elements \eqref{equation_hecke_gamma_b}, \eqref{equation_hecke_gamma_h} and \eqref{equation_hecke_gamma_0} gives the summation:
\begin{align*}
&p^{-2}\sum_{b'=0}^{p^2-1} \sum_{n\in\mathbb{Z}} a_{p^2c}^\pm(n) e\left(\frac{n(\tau+b'(p^2u_p^2-1))}{d_ph_cp^2}\right)e\left(\frac{-BR_c}{24}\right)\\
&=p^{-2}\sum_{n\in\mathbb{Z}} a_{p^2c}^\pm(n)q_{c}^{n/(d_pp^2)} \sum_{b'=0}^{p^2-1} e\left(  \frac{b'(p^2u_p^2-1)(n-p^2d_p\rho_c)}{d_ph_c p^2} \right).
\end{align*}
Since $pu_p\equiv1\pmod{d_pN}$ and $h_c\mid N$, we have $d_ph_c\mid p^2u_p^2-1$. Moreover, $p\nmid d_ph_c$ and $d_p\rho_c\in\mathbb Z$, so the inner exponential sum evaluates to $p^2$ when $p^2\mid n$ and $0$ otherwise.  Making the substitution $n = p^2 m$, the terms survive as $a_{p^2c}^\pm(p^2 m) q_c^{m/d_p}$, thereby giving the $a_{p^2c}^\pm(p^2 n)$ term.

We next consider the $p-1$ terms involving the expansion at the cusp $1/c$.
These include the $p-1$ elements $\alpha_b = \left(\begin{smallmatrix} 1 & b \\ 0 & p^2 \end{smallmatrix}\right)$ for $0\leq b\leq p^2-1$ such that $1+bc \equiv 0 \pmod p$ and $1+bc \not\equiv 0 \pmod{p^2}$.
For each such $b$, let $h\in\{1,\ldots,p-1\}$ satisfy
\[
(1 + bc)h \equiv -bp \pmod{p^2}.
\]
Write $H=h(p^2u_p^2-1)$ and consider the corresponding matrix
\[
\xi_b = \begin{pmatrix} \frac{p+(1+bc)cH}{p^2} & \frac{bp-(1+bc)H}{p^2} \\ c^2H  & p - cH \end{pmatrix}.
\]
The congruence implies $\xi_b \in \Gamma_0(N)$. Direct matrix multiplication verifies 
\[
\alpha_b\sigma_c=\xi_b\sigma_c\beta_H,
\]
which means
\[
(f|_{\kappa/2}\alpha_b)|_{\kappa/2}\sigma_c=\nu_g(\xi_b)(f|_{\kappa/2}\sigma_{c})|_{\kappa/2}\beta_H.
\]
Since $d_p\delta\mid H=h(p^2u_p^2-1)$ for every $\delta\mid N$, by Lemma~\ref{lemma_etamul_xi_b}, we obtain
\[
\nu_g(\xi_b)=\left(\frac{2\Delta H}{p}\right)\varepsilon_p^{-\kappa}e\left(\frac{bR-pHR_c}{24}\right).
\]
Summing over these $p-1$ representatives, the displayed decomposition contributes $\nu_g(\xi_b)$ itself, and hence the factor $\varepsilon_p^{-\kappa}$. Recalling that $H = h(p^2u_p^2-1)$, we obtain:
\[
p^{\lambda-3/2} \sum_{n\in\mathbb{Z}} a_c^\pm(n) q_{c}^{n/d_p}  \varepsilon_p^{-\kappa} \left(\frac{2\Delta}{p}\right) \sum_{h=1}^{p-1} \left(\frac{H}{p}\right) e\left(\frac{H(n-p^2d_p\rho_c)}{d_ph_cp}\right).
\]
Writing $H=h(p^2u_p^2-1)$ and using $d_ph_c\mid p^2u_p^2-1$ and $p\nmid d_ph_c$, the standard quadratic Gauss sum evaluates the inner sum as $\left(\frac{d_ph_cn}{p}\right)\varepsilon_p p^{1/2}$. Substituting this back, we use $\kappa = 2\lambda + 1$ to reduce the roots of unity:
\[
\varepsilon_p^{-\kappa} \cdot \varepsilon_p = \varepsilon_p^{1-\kappa} = \varepsilon_p^{-2\lambda} = \left(\frac{-1}{p}\right)^{-\lambda} = \left(\frac{(-1)^\lambda}{p}\right).
\]
Considering the Legendre symbols $\left(\frac{(-1)^\lambda}{p}\right) \left(\frac{d_p h_cn}{p}\right) \left(\frac{2\Delta}{p}\right) = \left(\frac{(-1)^\lambda(d_p/8)\Delta h_cn}{p}\right)$ and the weight factors, we get the second term of the theorem:
\[
\sum_{n\in\mathbb{Z}} \left(\frac{(-1)^\lambda(d_p/8)\Delta h_cn}{p}\right) p^{\lambda - 1} a_c^\pm(n) q_{c}^{n/d_p}.
\]

The remaining coset representative of $T(p^2)$ is $\alpha_{b} = \left(\begin{smallmatrix} 1 & b \\ 0 & p^2 \end{smallmatrix}\right)$, where $b$ is the unique integer modulo $p^2$ satisfying $1 + b c \equiv 0 \pmod{p^2}$.
Let
$
\sigma_b=\begin{pmatrix}
(1+bc)/p^2-u_p^2bc & b\\
(1-p^2u_p^2)c & p^2
\end{pmatrix}.
$
We have
\[
\alpha_b\sigma_c=\sigma_b\sigma_{u_p^2c}\sigma,
\]
with cocycle factor $1$: as $\operatorname{Im}\tau\to\infty$, the expression inside the square root associated with $\sigma_b$ tends to $u_p^{-2}>0$.
By Proposition~\ref{proposition_eta_multiplier_p_square}, Proposition~\ref{proposition_etamul_squares} and $pu_p\equiv 1\pmod{d_pN}$,
\begin{align*}
\nu_g(\sigma_b)&=\nu_g(p^2,b,(1-p^2u_p^2)c,(1+bc)/p^2-u_p^2bc)\\
&=e\left(\frac{bR}{24}\right)e\left(\frac{(p^2u_p^2-1)c R_1}{24}\right)=e\left(\frac{bR}{24}\right).
\end{align*}
The slash action of $\sigma$ contributes $p^{\kappa/2}$ and replaces $\tau$ by $p^2\tau$. Combining this with the normalization $p^{\lambda-3/2}$ and reindexing the sum, we see that the third term is
\[
\sum_{n \in \mathbb{Z}} p^{2\lambda-1} a_{u_p^2c}^\pm(n/p^2) q_{c}^{n/d_p}.
\]
This proves the formulas for the holomorphic coefficients and for the non-holomorphic coefficients with $n\ne0$.

It remains to check the $n=0$ term in the non-holomorphic part. The first group of representatives replaces $y$ by $y/p^2$, so its coefficient acquires the factor $p^{\kappa-2}$. The last representative replaces $y$ by $p^2y$, and its original
factor $p^{\kappa-2}$ is therefore cancelled. The middle Gauss sum vanishes for $n=0$. This gives
\[
A_c^-(0)=p^{\kappa-2}a_{p^2c}^-(0)+a_{u_p^2c}^-(0),
\]
as asserted.
\end{proof}

\begin{proof}[Proof of Lemma~\ref{lemma_etamul_gamma_b}]
Let 
\[
m =\frac{\gcd(c,\delta)^2B}{\delta}.
\]
We first consider $b=0$. In this case, the matrix is
\[
\gamma_{0,\delta} = \begin{pmatrix} 1+Bc & -B\delta \\ Bc^2/\delta & 1-Bc \end{pmatrix}.
\]
Substituting $Bc = m c_1 \delta_1$, $c = \gcd(c,\delta) c_1$, and $\delta = \gcd(c,\delta) \delta_1$, the matrix reduces to:
\[
\gamma_{0,\delta} = \begin{pmatrix} 1 + m c_1 \delta_1 & -m\delta_1^2 \\ m c_1^2 & 1 - m c_1 \delta_1 \end{pmatrix}.
\]
Because $\gcd(c_1, \delta_1) = 1$, by Proposition~\ref{proposition_easy_pro_eta} (8), we have
\[
\nu_{\eta}(\gamma_{0,\delta})=\nu_{\eta}(1 + m c_1\delta_1,-m\delta_1^2,m c_1^2, 1 - m c_1 \delta_1)=e(-m/24).
\]
To generalize to $b \neq 0$, we note that $\gamma_{b,\delta}$ is a left-translation of $\gamma_{0,\delta}$ by $b\delta$. By Proposition~\ref{proposition_easy_pro_eta} (6),
\[
\nu_{\eta}(\gamma_{b,\delta})=e(b\delta/24)e(-m/24)=e\left( \frac{b\delta}{24} - \frac{B\gcd(c,\delta)^2}{24\delta} \right).
\]
This completes the proof.
\end{proof}

\begin{proof}[Proof of Lemma~\ref{lemma_etamul_gamma_h}]
The congruence hypotheses ensure that all entries are integral, and direct calculation shows that the matrix belongs to $\SL_2(\mathbb{Z})$. Write $B_0=B/\delta$ and
\[
\gamma_{h,\delta}=\begin{pmatrix}A&J\\C&D\end{pmatrix}, \qquad A=p+pBc+Bc^2h,\quad C=pB_0c^2.
\]
Since $24\mid B_0$, we have
\[
A\equiv p\pmod{24},\qquad
J\equiv h\delta\pmod{24},\qquad
C\equiv0\pmod{24}.
\]
Also $A,C>0$. Proposition~\ref{proposition_etamul_squares} yields
\begin{equation}\label{equation_gamma_h_short_start}
\nu_\eta(\gamma_{h,\delta})=\left(\frac{2C}{A}\right)\varepsilon_p^{-1}
e\left(\frac{ph\delta}{24}\right).
\end{equation}

We evaluate the symbol in this expression. Since $\gcd(A,c)=1$,
\[
\left(\frac{2C}{A}\right)=\left(\frac{2}{A}\right)
\left(\frac{p}{A}\right)\left(\frac{B_0}{A}\right).
\]
The congruences $A\equiv p\pmod8$ and $A\equiv p\pmod{B_0}$, together with quadratic reciprocity, imply
\[
\left(\frac{2}{A}\right)=\left(\frac{2}{p}\right),\qquad
\left(\frac{B_0}{A}\right)=\left(\frac{B_0}{p}\right).
\]
For the second equality, one may write $B_0=2^v u$ with $u$ odd,
use $\left(\frac{A}{u}\right)=\left(\frac{p}{u}\right)$, and note that $A\equiv p\pmod4$ makes
the reciprocity signs identical. Moreover,
\[
\left(\frac{p}{A}\right)=\left(\frac{A}{p}\right)
(-1)^{\frac{p-1}{2}\frac{A-1}{2}}=\left(\frac{hc}{p}\right)\left(\frac{-1}{p}\right),
\]
because $A\equiv hc\pmod p$ and $A\equiv p\pmod4$.
Finally, $B_0\delta c\equiv1\pmod p$ implies
$\left(\frac{B_0}{p}\right)=\left(\frac{\delta c}{p}\right)$. Hence
\[
\left(\frac{2C}{A}\right)=\left(\frac{-1}{p}\right)\left(\frac{2h\delta}{p}\right).
\]
Substituting this into \eqref{equation_gamma_h_short_start} and using
$\left(\frac{-1}{p}\right)=\varepsilon_p^2$ gives the asserted formula.
For $p=3$ and $8\mid B_0$, the three entry congruences hold modulo $8$.
By Lemma~\ref{lemma_eta_cubes_three}, we obtain
\[
\nu_\eta(\gamma_{h,\delta})^3=\left(\frac{2C}{A}\right)\varepsilon_p^{-3}e(ph\delta/8).
\]
The symbol calculation above uses only congruences modulo $8$, $4$, $B_0$, and $p$, so remains valid. Since $\varepsilon_p^2\varepsilon_p^{-3}=\varepsilon_p^3$, this is exactly the cube of the stated formula.
\end{proof}

\begin{proof}[Proof of Lemma~\ref{lemma_etamul_xi_b}]
Write $H_0=H/\delta$ and $\xi_{b,\delta}=\left(\begin{smallmatrix}A&B\delta\\C&D\end{smallmatrix}\right)$, where $C=c^2H_0$. The congruence hypotheses imply that all entries are integral, and direct calculation shows that it belongs to $\SL_2(\mathbb{Z})$. If $k=(1+bc)/p$, then $p\nmid k$ and $kH\equiv b\pmod p$. In particular, $p\nmid H$. Moreover,
\[
pA=1+kcH,
\]
so $pA\equiv1\pmod{H_0}$.

Since $24\mid H_0$ and $p^2\equiv1\pmod{24}$, we have
\[
A\equiv p\pmod{24},\qquad
B\delta\equiv bp\delta\pmod{24},\qquad
C\equiv0\pmod{24}.
\]
Also $A>0$ and $C>0$. Therefore by Proposition~\ref{proposition_etamul_squares},
we obtain
\[
\nu_\eta(\xi_{b,\delta})=\left(\frac{2C}{A}\right)\varepsilon_A^{-1}
e\left(\frac{A(B\delta-C)}{24}\right)=\left(\frac{2H_0}{A}\right)\varepsilon_p^{-1}e\left(\frac{b\delta}{24}\right),
\]
where we used $\gcd(A,c)=1$ to remove $c^2$ from the symbol.

Write $H_0=2^v u$ with $u$ odd. From $A\equiv p\pmod8$, we obtain $\left(\frac{2}{A}\right)=\left(\frac{2}{p}\right)$. From $pA\equiv1\pmod u$, we obtain $\left(\frac{A}{u}\right)=\left(\frac{p}{u}\right)$. Quadratic reciprocity, together with $A\equiv p\pmod4$, now gives $\left(\frac{u}{A}\right)=\left(\frac{u}{p}\right)$. Consequently,
\[
\left(\frac{2H_0}{A}\right)=\left(\frac{2H_0}{p}\right)=\left(\frac{2\delta H}{p}\right).
\]
Finally, $24\mid H_0$ implies
\[
e\left(-\frac{pH\gcd(c,\delta)^2}{24\delta}\right)=1,
\]
which yields the stated formula. For $p=3$ and $8\mid H_0$, the entry congruences instead hold modulo $8$. By Lemma~\ref{lemma_eta_cubes_three}, we know
\[
\nu_\eta(\xi_{b,\delta})^3=\left(\frac{2H_0}{A}\right)\varepsilon_p^{-3}e(b\delta/8).
\]
The reciprocity calculation is unchanged, and the omitted exponential factor has cube $1$ because $8\delta\mid H$. This proves the cubed version.
\end{proof}

\section{Proofs of the main results}\label{section_main_proofs}
\subsection{Preliminaries}
Let $f(\tau) \in H_{\kappa/2}^!(\Gamma_0(N), \overline{\nu}_g)$, where $\overline{\nu}_g$ is associated with the eta-quotient $g(\tau) = \prod_{\delta\mid N} \eta(\delta\tau)^{r_\delta}$. We assume \eqref{equation_prime_hypotheses} and use the notation in \eqref{equation_main_arithmetic_notation} and \eqref{equation_uniform_local_indices}. Here $d\rho_c=dh_cR_c/24\in\mathbb Z$. Choose $u_\ell>0$ with $\ell u_\ell\equiv1\pmod{dN}$. For the conjugate multiplier we apply Section~\ref{section_hecke_operators} to $g^{-1}$: the parameter $R_c$ changes sign, while $\Delta^{-1}$ and $\Delta$ give the same Legendre symbols. The compatibility condition is $\kappa\equiv-r\pmod4$.  We set $T_1=I$ and $T_{\ell^2}=T(\ell^2)$. For $m\geq2$, we define the operators recursively by
\begin{equation}\label{equation_hecke_recurrence}
T_{\ell^{2m}} = T_{\ell^2}T_{\ell^{2m-2}} - \ell^{2\lambda-1} T_{\ell^{2m-4}}.
\end{equation}

\begin{lemma}\label{lemma_twisted_hecke_principal}
Fix $\chi(\ell)\in\{1,-1\}$. Suppose that, at every cusp,
\[
n<0,\quad a_c(n)\ne0\quad\Longrightarrow\quad
\epsilon_\ell(-h_cn)=\chi(\ell).
\]
Then, for every $m\ge1$,
\begin{equation}\label{equation_hecke_principal_difference}
\begin{aligned}
&\mathcal P_c^-(T_{\ell^{2m}}f-\ell^{\lambda-1}\chi(\ell)T_{\ell^{2m-2}}f)\\
&\qquad=\ell^{(2\lambda-1)m}\sum_{n<0}
a_{u_\ell^{2m}c}(n)q_c^{\ell^{2m}n/d}.
\end{aligned}
\end{equation}
\end{lemma}

\begin{proof}
Every nonzero initial polar coefficient has $\ell\nmid n$. The quantities $h_c$ and $\rho_c$ are unchanged under $c\mapsto\ell^2c$ and $c\mapsto u_\ell^2c$. The coefficient formula in Theorem~\ref{theorem_hecke_operator} does not mix negative and nonnegative indices. We claim that, for $m\ge0$,
\begin{equation}\label{equation_hecke_principal}
\mathcal P_c^-(T_{\ell^{2m}}f)=
\sum_{n<0}\sum_{j=0}^m
(\ell^{\lambda-1}\chi(\ell))^{m-j}\ell^{(2\lambda-1)j}
a_{u_\ell^{2j}c}(n)q_c^{\ell^{2j}n/d}.
\end{equation}

The formula is immediate for $m=0$. For $m=1$, the term $a_{\ell^2c}(\ell^2n)$ contributes nothing to the principal part, the middle term multiplies each initial polar coefficient by $\ell^{\lambda-1}\chi(\ell)$, and the last term multiplies its exponent by $\ell^2$ and its coefficient by $\ell^{2\lambda-1}$. This proves \eqref{equation_hecke_principal} for $m=1$.

For the induction step, apply $T_{\ell^2}$ to the expression for $m-1$. The terms that divide the exponent by $\ell^2$ have $j\ge1$. By the periodicity $f_{c+dN}=f_c$ established in the proof of Theorem~\ref{theorem_hecke_operator} and by $\ell^2u_\ell^2\equiv1\pmod{dN}$,
\[
a_{u_\ell^{2j}\ell^2c}(n)=a_{u_\ell^{2j-2}c}(n).
\]
Thus these terms are exactly $\ell^{2\lambda-1}$ times the expression for $m-2$, and cancel in \eqref{equation_hecke_recurrence}. The middle term is zero for $j\ge1$, since its quadratic symbol has argument divisible by $\ell$. At $j=0$ it supplies the $j=0$ term for $m$, and the terms multiplying the exponent by $\ell^2$ supply the terms with $1\le j\le m$. This proves \eqref{equation_hecke_principal}. Subtracting $\ell^{\lambda-1}\chi(\ell)$ times the formula for $m-1$ leaves only the terms with $j=m$, proving \eqref{equation_hecke_principal_difference}.
\end{proof}

We use the following variant of \cite[Proposition~2.3]{zbMATH07595066}, in which the third identity also covers the case $\ell^2\mid n$.
\begin{lemma} \label{lemma_fourier_coeffs_fm_infty}
Fix $\chi(\ell)\in\{1,-1\}$.
Suppose that the Fourier expansion of the holomorphic part of $f$ at $\infty$ is given by
\[
f^{+}(\tau) = \sum_{n \gg -\infty} a(n) q^{n/d}.
\]
Define $c_m(n)$ for $m\ge1$ by
\[
\left(T_{\ell^{2m}}f(\tau)-\ell^{\lambda-1}\chi(\ell)T_{\ell^{2m-2}}f\right)^{+} = \sum_{n \gg -\infty} c_m(n) q^{n/d}.
\]
Put $c_0(n)=a(n)$. For $m\ge1$, the Fourier coefficients $c_m(n)$ satisfy the following identities, with the convention that any coefficient with a non-integer index is zero:
\begin{enumerate}
\item For any $n$, we have
\[
c_m(\ell^2 n) - \ell^{2\lambda-1}c_{m-1}(n) =a(\ell^{2m+2} n) - \ell^{\lambda-1}\chi(\ell) a(\ell^{2m} n).
\]

\item If $\ell \nmid n$, then
\begin{align*}
c_m(n) &= a(\ell^{2m} n) \\
&+ \left[ 1 - \chi(\ell)\epsilon_\ell(-n) \right] \sum_{j=1}^m (-1)^j
\left(\ell^{\lambda-1} \chi(\ell)\right)^j a(\ell^{2m-2j} n).
\end{align*}

\item If $\ell \mid n$, then
\[
c_m(n) = a(\ell^{2m} n) - \ell^{\lambda-1}\chi(\ell) a(\ell^{2m-2} n) + \ell^{2\lambda-1} c_{m-1}\left( \frac{n}{\ell^2} \right).
\]
\end{enumerate}
\end{lemma}

\begin{proof}
Write $\mathcal F_0=f$ and
$\mathcal F_m=T_{\ell^{2m}}f-\ell^{\lambda-1}\chi(\ell)T_{\ell^{2m-2}}f$ for $m\ge1$.
For $F(\tau)=\sum C(n)q^{n/d}$, let $(T_{\ell^2}F)(n)$ denote the coefficient of $q^{n/d}$ in $T_{\ell^2}F$. By Theorem~\ref{theorem_hecke_operator}, we have
\[
(T_{\ell^2} F)(n) = C(\ell^2 n) + \ell^{\lambda-1} \epsilon_\ell(-n) C(n) + \ell^{2\lambda-1} C\left( \frac{n}{\ell^2} \right).
\]
In particular, since the Legendre symbol vanishes, we have
\begin{equation}\label{equation_hecke_operator_ell_square_index}
(T_{\ell^2} F)(\ell^2 n) = C(\ell^4 n) + \ell^{2\lambda-1} C(n).
\end{equation}
For Part (1), we proceed by induction on $m$. Applying \eqref{equation_hecke_operator_ell_square_index} to $f^{+}$ yields $c_1(\ell^2 n) = a(\ell^4 n) + \ell^{2\lambda-1} a(n) - \ell^{\lambda-1} \chi(\ell) a(\ell^2 n)$. Subtracting $\ell^{2\lambda-1} c_0(n) = \ell^{2\lambda-1} a(n)$ verifies the base case $m=1$.

Assuming the identity for $m-1$, the Hecke relation $\mathcal F_m =T_{\ell^2} \mathcal F_{m-1} - \ell^{2\lambda-1} \mathcal F_{m-2}$ evaluated at $\ell^2 n$ directly shows
\[
c_m(\ell^2 n) = c_{m-1}(\ell^4 n)+ \ell^{2\lambda-1} c_{m-1}(n) - \ell^{2\lambda-1} c_{m-2}(\ell^2 n).
\]
Rearranging this equation and applying the induction hypothesis gives the formula for $m$, completing Part (1). 

For Parts (2) and (3), the case $m=1$ is directly from the Hecke operator formula.  For $m\ge2$, apply it to $\mathcal F_{m-1}$ and use Part (1) to obtain
\begin{align*}
 c_m(n)
 &=a(\ell^{2m}n)
   -\ell^{\lambda-1}\chi(\ell)a(\ell^{2m-2}n)\\
 &\quad+\ell^{\lambda-1}\epsilon_\ell(-n) c_{m-1}(n)
   +\ell^{2\lambda-1}c_{m-1}\left(\frac{n}{\ell^2}\right).
\end{align*}
This identity also holds for $m=1$ with $c_0(n)=a(n)$. If $\ell\nmid n$, the last index is nonintegral. Iteration of the remaining recurrence proves Part (2). If $\ell\mid n$, then $\epsilon_\ell(-n)=0$, giving Part (3).
\end{proof}

\subsection{A vanishing criterion}
\label{subsection_dimension_vanishing}
We recall the dimension formula for the space of cusp forms. For even $k\ge2$, the standard formula \cite[Chapter~3]{diamond2005first} is
\begin{equation}\label{equation_cusp_dimension_formula}
\begin{aligned}
d_k(N)&=\frac{k-1}{12}\mu_0(N)-\frac12\varepsilon_\infty(N)+\left(\frac14-\left\{\frac k4\right\}\right)\varepsilon_2(N)\\
&\quad+\left(\frac13-\left\{\frac k3\right\}\right)\varepsilon_3(N)+\delta_{k,2},
\end{aligned}
\end{equation}
where $\{x\}=x-\lfloor x\rfloor$, and $\mu_0(N):=\left[\mathrm{SL}_2(\mathbb{Z}):\Gamma_0(N)\right]$ denotes the index of $\Gamma_0(N)$ in $\mathrm{SL}_2(\mathbb{Z})$; and $\varepsilon_2(N)$, $\varepsilon_3(N)$ count elliptic points of orders $2$ and $3$, respectively; and $\varepsilon_\infty(N)$ counts the number of cusps. Combining \eqref{equation_cusp_dimension_formula} with its weight $2$ specialization, we obtain
\begin{equation}\label{equation_dimension_bound_identity}
\begin{aligned}
&\frac{k\mu_0(N)}{12} -\left\{\frac k4\right\}\varepsilon_2(N)-\left\{\frac k3\right\}\varepsilon_3(N)-\varepsilon_\infty(N)+1\\
&=d_k(N)+d_2(N)-\delta_{k,2}.
\end{aligned}
\end{equation}
For a form in $M_k(\Gamma_0(N))$ with nonzero reduction, the vanishing-order bound of Ahlgren, Masri, and Rouse \cite[Theorem~4.1]{ahlgren2009vanishing} is
$d_k(N)+d_2(N)+\varepsilon_\infty(N)-\delta_{k,2}-1$. For cusp forms it is $d_k(N)+d_2(N)-\delta_{k,2}$. Following their argument, we retain the vanishing at all cusps and allow weakly holomorphic forms in the criterion below.

\begin{lemma}\label{lemma_integral_cusp_vanishing}
Let $k\ge2$ be even, let $m\ge1$, and let $\ell\ge3$ be prime with $\ell\nmid N$. Let $H\in M_k^!(\Gamma_0(N))$, and use the local coefficient ring $\mathcal O$ fixed in Section~\ref{subsection_cusp_expansions}. Suppose that every cusp expansion
\[
H|_k\sigma_{\mathfrak a}=\sum_{n\gg-\infty}h_{\mathfrak a}(n)q_{\mathfrak a}^n
\]
has coefficients in $\mathcal O$. At each cusp choose an integer $b_{\mathfrak a}\in\mathbb Z$ such that
\[
h_{\mathfrak a}(n)\in\ell^m\mathcal O\quad\text{for every }n<b_{\mathfrak a}.
\]
If
\begin{equation}\label{equation_all_cusp_vanishing_condition}
\sum_{\mathfrak a}(b_{\mathfrak a}-1)\ge d_k(N)+d_2(N)-\delta_{k,2},
\end{equation}
then $H\equiv0\pmod{\ell^m}$ coefficientwise at every cusp.
\end{lemma}

\begin{proof}
Enlarge $\mathcal O$ by a finite extension so that the cusp expansions and the required roots of unity are defined over it. Such an extension is faithfully flat and preserves the divisibility assertion. Let $\varpi$ be a uniformizer, write $\ell\mathcal O=\varpi^e\mathcal O$, and denote reduction modulo $\varpi$ by a bar.

We first prove the assertion modulo $\varpi$. Let $F\in M_k^!(\Gamma_0(N))$ have Fourier coefficients in $\mathcal O$ at all cusps. At each cusp $\mathfrak a$, assume that the coefficient of $q_{\mathfrak a}^n$ belongs to $\varpi\mathcal O$ for every $n<b_{\mathfrak a}$. In this proof, $\Delta(\tau)=\eta(\tau)^{24}$ denotes the Ramanujan Delta function. Choose $a\ge0$ sufficiently large that
\[
G=F\Delta^a\in M_{k+12a}(\Gamma_0(N)),\qquad
ah_{\mathfrak a}+b_{\mathfrak a}\ge1
\quad\text{at every cusp }\mathfrak a.
\]
The integral $q$-expansion principle at good reduction, using $\ell\nmid N$, implies that $G$ is defined over $\mathcal O$ and that, if one cusp expansion has nonzero reduction, every cusp expansion has nonzero reduction; see \cite[\S1.6]{katz1973p}. Suppose that this is the case. Since $\Delta|_{12}\sigma_{\mathfrak a}$ begins with $q_{\mathfrak a}^{h_{\mathfrak a}}$ with coefficient $1$, we have
\[
\operatorname{ord}_{q_{\mathfrak a}}
\overline{G|_{k+12a}\sigma_{\mathfrak a}}\ge ah_{\mathfrak a}+b_{\mathfrak a}.
\]

Write $w=k+12a$. Following the proof of \cite[Theorem~4.1]{ahlgren2009vanishing}, we define the norm
\[
P=\prod_{\gamma\in\Gamma_0(N)\backslash\mathrm{SL}_2(\mathbb Z)}
G|_w\gamma\in M_{w\mu_0(N)}(\mathrm{SL}_2(\mathbb Z)).
\]
Every factor has coefficients in $\mathcal O$ and nonzero reduction. A cusp of width $h_{\mathfrak a}$ contributes $h_{\mathfrak a}$ translated factors, so
\[
\operatorname{ord}_q\overline P=\sum_{\mathfrak a}\operatorname{ord}_{q_{\mathfrak a}}\overline{G|_w\sigma_{\mathfrak a}}.
\]
Let
\[
A_2=2\left\{\frac w4\right\}\varepsilon_2(N),\qquad
A_3=3\left\{\frac w3\right\}\varepsilon_3(N).
\]
At an elliptic point of order $2$ or $3$, the transformation law forces the vanishing order to be at least $2\{w/4\}$ or $3\{w/3\}$, respectively. Thus
\[
Q=\frac{P}{E_4^{A_3}E_6^{A_2}}
\]
is a holomorphic modular form of level one: the denominator removes
only forced elliptic zeros in characteristic zero. Since
$E_4,E_6\in1+q\mathbb Z[[q]]$, their formal inverses have integral
coefficients. Thus $Q$ has coefficients in $\mathcal O$ and nonzero reduction.
This argument also applies when $\ell=3$. Applying Sturm's theorem at level one \cite{sturm1987congruence} to $Q$, we have
\begin{equation}\label{equation_all_cusp_order_bound}
\sum_{\mathfrak a}\operatorname{ord}_{q_{\mathfrak a}}
\overline{G|_w\sigma_{\mathfrak a}}\le \frac{w\mu_0(N)}{12}-\left\{\frac w4\right\}\varepsilon_2(N)-\left\{\frac w3\right\}\varepsilon_3(N).
\end{equation}

Since $\sum_{\mathfrak a}h_{\mathfrak a}=\mu_0(N)$ and
$w=k+12a$, the fractional parts in this bound agree with those for
$k$. Subtracting $a\mu_0(N)$ and using
\eqref{equation_dimension_bound_identity}, we obtain
\[
\sum_{\mathfrak a}(b_{\mathfrak a}-1)\le d_k(N)+d_2(N)-\delta_{k,2}-1,
\]
contrary to \eqref{equation_all_cusp_vanishing_condition}. Thus every coefficient in every cusp expansion of $F$ is divisible by $\varpi$.

Apply this assertion successively to $H/\varpi^j$ for $0\le j<em$. At each step, the induction shows that the Fourier coefficients at every cusp belong to $\mathcal O$, while the original hypothesis ensures that the coefficient of $q_{\mathfrak a}^n$ is divisible by $\varpi$ whenever $n<b_{\mathfrak a}$. Thus every coefficient at every cusp belongs to $\varpi^{em}\mathcal O=\ell^m\mathcal O$.
\end{proof}

\begin{remark}\label{remark_computable_cusp_criterion}
The definition \eqref{equation_signed_eta_cusp_excess} can also be written as
\[
\mathcal W_N(g)=\sum_{c\mid N}\varphi(\gcd(c,N/c))
\bigl(\lfloor\rho_c\rfloor-\mathbf1_{\rho_c\in\mathbb Z_{>0}}\bigr),
\]
where $\mathbf1$ is the indicator function. Indeed, by the  classification of cusp \cite[Chapter~3]{diamond2005first}, there are $\varphi(\gcd(c,N/c))$ cusps with denominator $c$, and the eta-quotient has the same order at each of them. The local bounds used below are
\[
b_{\mathfrak a}=\begin{cases}
\lceil\rho_{\mathfrak a}\rceil,&\rho_{\mathfrak a}>0,\\
1+\lfloor\rho_{\mathfrak a}\rfloor,&\rho_{\mathfrak a}\le0.
\end{cases}
\]
Since $b_{\mathfrak a}-1=\lfloor\rho_{\mathfrak a}\rfloor-\mathbf1_{\rho_{\mathfrak a}\in\mathbb Z_{>0}}$, their sum is $\mathcal W_N(g)$. At a positive integer order, the subtraction of $1$ accounts for the possible constant term of the Hecke difference. For holomorphic $g$ these bounds are $b_{\mathfrak a}=\max\{1,\lceil\rho_{\mathfrak a}\rceil\}$.
\end{remark}

\begin{proposition}\label{proposition_cusp_dimension_equivalence}
Suppose that $\kappa\ge5$ is odd and that $k=(\kappa+r)/2$ is even with $k\ge2$. Then
\begin{equation}\label{equation_half_integral_cusp_dimension}
\dim S_{\kappa/2}(\Gamma_0(N),\overline\nu_g)
=d_k(N)-\delta_{k,2}-\sum_{\mathfrak a}\lfloor\rho_{\mathfrak a}\rfloor.
\end{equation}
Consequently, \eqref{equation_dimension_cusp_criterion} is equivalent to
\[
d_2(N)=0,\qquad S_{\kappa/2}(\Gamma_0(N),\overline\nu_g)=\{0\},
\qquad \rho_{\mathfrak a}\notin\mathbb Z_{>0}\quad\text{at every cusp}.
\]
\end{proposition}

\begin{proof}
Multiplication by $g$ identifies $S_{\kappa/2}(\Gamma_0(N),\overline\nu_g)$ with the space of weight-$k$ forms holomorphic on $\mathbb H$ and of order at least $1+\lfloor\rho_{\mathfrak a}\rfloor$ at each cusp. This is the space of sections of a line bundle $\mathcal L$ on $X_0(N)$; see the divisor description in \cite[Lemma~4.1]{zhu2024dimension}. Removing the forced elliptic zeros and using \eqref{equation_dimension_bound_identity}, we obtain
\[
\begin{aligned}
\deg\mathcal L&=\frac{k\mu_0(N)}{12}-\left\{\frac k4\right\}\varepsilon_2(N)-\left\{\frac k3\right\}\varepsilon_3(N)-\sum_{\mathfrak a}(1+\lfloor\rho_{\mathfrak a}\rfloor)\\
&=d_k(N)+d_2(N)-\delta_{k,2}-1-\sum_{\mathfrak a}\lfloor\rho_{\mathfrak a}\rfloor.
\end{aligned}
\]

Under the pairing $H,u\mapsto Hu\,d\tau$, sections of $K_{X_0(N)}\otimes\mathcal L^{-1}$ correspond to forms $u$ of weight $2-k$, holomorphic on $\mathbb H$, with $\operatorname{ord}_{\mathfrak a}(u)\ge-\lfloor\rho_{\mathfrak a}\rfloor$. Indeed, at an elliptic point of order $e$, the complementary local vanishing orders are $a$ and $e-1-a$, with $0\le a<e$, so their product times $d\tau$ is a regular differential in the local coordinate $z^e$. At a cusp, $d\tau$ is a nonzero multiple of $dq_{\mathfrak a}/q_{\mathfrak a}$, which gives the stated bound. Multiplication by $g$ identifies these forms $u$ with $M_{2-\kappa/2}(\Gamma_0(N),\nu_g)$: the integral cusp exponents of $u$ make the cusp bound equivalent to holomorphy of $gu$, and $g$ has no zeros on $\mathbb H$. This space is zero, since its weight is negative and its multiplier has finite order. By the Riemann--Roch theorem, we obtain
\[
\dim S_{\kappa/2}(\Gamma_0(N),\overline\nu_g)=\deg\mathcal L+1-d_2(N),
\]
and thus prove \eqref{equation_half_integral_cusp_dimension}.
Substituting into \eqref{equation_dimension_cusp_criterion}, we obtain
\[
\dim S_{\kappa/2}(\Gamma_0(N),\overline\nu_g)+d_2(N)+\#\{\mathfrak a:\rho_{\mathfrak a}\in\mathbb Z_{>0}\}\le0.
\]
\end{proof}

\subsection{Weakly holomorphic forms}
\begin{proof}[Proof of Theorem~\ref{theorem_weakly_holomorphic}]
Define $F_m=T_{\ell^{2m}}f-\ell^{\lambda-1}\chi(\ell) T_{\ell^{2m-2}}f$
and $b=(\kappa-2)m$. By Lemma~\ref{lemma_cc_magnitude}, $f$ has $\ell$-integral coefficients at every cusp. The coefficient formula in Theorem~\ref{theorem_hecke_operator} preserves $\ell$-integrality because $\kappa\ge3$. Thus the Fourier coefficients of $F_m$ and $H_m=gF_m\in M_k^!(\Gamma_0(N))$ at every cusp belong to $\mathcal O$. By (\ref{hypothesis_weak_principal}) and Lemma~\ref{lemma_twisted_hecke_principal}, all coefficients of the principal part of $F_m$ at every cusp are divisible by $\ell^b$.

At a cusp with $\rho_c\le0$, the constant term of $f_c$ is zero: this follows from (\ref{hypothesis_weak_constant}) when $\rho_c$ is integral and from the translation law otherwise. The same holds throughout its Hecke orbit. The coefficient formula at index zero involves only constant terms, so $F_m$ also has zero constant term there. Modulo $\ell^b$, its expansion therefore has strictly positive exponents in $\mathbb Z-\rho_c$. Consequently the coefficient of $q_c^n$ in $H_m|_k\sigma_c$ belongs to $\ell^b\mathcal O$ for every $n<1+\lfloor\rho_c\rfloor$. At a cusp with $\rho_c>0$, the vanishing of the principal part modulo $\ell^b$ gives the same divisibility for every $n<\lceil\rho_c\rceil$.

Take $b_{\mathfrak a}$ as in Remark \ref{remark_computable_cusp_criterion}. Then
$\sum_{\mathfrak a}(b_{\mathfrak a}-1)=\mathcal W_N(g)$, so $H_m\equiv0\pmod{\ell^b}$ by
\eqref{equation_dimension_cusp_criterion} and Lemma
\ref{lemma_integral_cusp_vanishing}. By Lemma~\ref{lemma_cc_magnitude}, division by $g_c$ preserves coefficientwise divisibility. Thus
\begin{equation}\label{equation_weak_hecke_congruence}
T_{\ell^{2m}}f-\ell^{\lambda-1}\chi(\ell) T_{\ell^{2m-2}}f
\equiv0\pmod{\ell^{(\kappa-2)m}}
\end{equation}
at every cusp.

Write $F_m=\sum_n c_m(n)q^{n/d}$ at infinity and set $c_0(n)=a(n)$. If $\chi(\ell)\epsilon_\ell(-n)=1$, Part (2) of Lemma~\ref{lemma_fourier_coeffs_fm_infty} implies
\[
a(\ell^{2m}n)=c_m(n)\equiv0\pmod{\ell^{(\kappa-2)m}},
\]
which proves the first congruence. For part (\ref{conclusion_weak_recurrence}), put $m=\lceil j/2\rceil$. If $j=2m$, Part (1) of the lemma expresses the difference in \eqref{equation_main_unified_recurrence} as
\[
c_m(\ell^2n)-\ell^{\kappa-2}c_{m-1}(n).
\]
If $j=2m-1$, Part (3), applied at $\ell n$, expresses the same difference as
\[
c_m(\ell n)-\ell^{\kappa-2}c_{m-1}(n/\ell).
\]
Both expressions are divisible by $\ell^{(\kappa-2)m}$, using \eqref{equation_weak_hecke_congruence} for $m$ and $m-1$. For $m=1$, use $c_0(n)=a(n)\in\mathcal O$; coefficients with nonintegral indices are zero. This proves the recurrence.
\end{proof}

\begin{remark}\label{remark_weak_quotient_criterion}
Suppose $f=E/g$, where $E\in M_k(\Gamma_0(N))$ has $\ell$-integral coefficients at infinity. Write $g_c=C_cq_c^{\rho_c}U_c(q_c)$ as in Lemma~\ref{lemma_cc_magnitude}. At a cusp with $\rho_c<0$, the quotient vanishes. At a cusp with $\rho_c=0$, (\ref{hypothesis_weak_constant}) is equivalent to $E_c(0)=0$. If $\rho_c>0$, check (\ref{hypothesis_weak_principal}) on the finitely many negative powers in $q_c^{-\rho_c}E_c/(C_cU_c)$. When $0<\rho_c\le1$, the only possible polar coefficient is $a_c(-d\rho_c)=E_c(0)/C_c$, so the condition reduces to
\[
E_c(0)\ne0\quad\Longrightarrow\quad
\epsilon_\ell(dh_c\rho_c)=\chi(\ell).
\]
The weight restriction, \eqref{equation_prime_hypotheses}, and \eqref{equation_dimension_cusp_criterion} are still required.
\end{remark}

\subsection{Harmonic weak Maass forms of weight \texorpdfstring{$3/2$}{3/2}}\label{subsection_harmonic_weak_maass}
We first establish the two facts about the shadow that will be used in the proof.

\begin{lemma}\label{lemma_eta_hecke_eigenform}
Let $g\in M_{1/2}(\Gamma_0(N),\nu_g)$ be a holomorphic eta-quotient,
and assume \eqref{equation_prime_hypotheses} with $\ell=p$.
Then $g$ is a Hecke eigenform, with
\begin{equation}\label{equation_eta_hecke_eigenvalue}
T(p^2)g=(1+p^{-1})\chi(p)g,\qquad \chi(p)\in\{1,-1\}.
\end{equation}
\end{lemma}
\begin{proof}
By Mersmann's classification, every holomorphic eta-quotient of weight $1/2$ is an integral rescaling of one of fourteen primitive eta-quotients; see \cite[Theorem~1]{bhattacharya2017holomorphic}. Here an eta-quotient is primitive if it is not of the form $t(d\tau)$ with $d>1$ and another eta-quotient $t$. For the corresponding theta expansions, see \cite[\S3]{bhattacharya2017holomorphic} and \cite[Theorems~1.1 and 1.2]{lemkeoliver2013eta}.

For a primitive quotient, write $g=\sum_{n\ge0}a(n)q^{n/d_p}$, with $a(x)=0$ for nonintegral $x$. The cited theta expansions are supported on a single square class. Multiplication of the theta summation index by $p$ gives a sign $\chi(p)$ for which
\[
\begin{aligned}
a(p^2n)&=\chi(p)a(n),\\
\left(\frac{(d_p/8)\Delta n}{p}\right)a(n)&=\chi(p)a(n)
&&(p\nmid n),\\
a(n/p^2)&=\chi(p)a(n)&&(p\mid n).
\end{aligned}
\]
For $p\ge5$, these identities follow directly from the periodic coefficients in those theta expansions \cite[Appendix~A, Table~1]{wang2026explicit}. For $p=3$, the restrictions $3\nmid N$ and $3\mid R$ leave only the primitive quotients
\[
\frac{\eta(\tau)^2}{\eta(2\tau)},\qquad
\frac{\eta(2\tau)^2}{\eta(\tau)},\qquad
\frac{\eta(2\tau)^5}{\eta(\tau)^2\eta(4\tau)^2},\qquad
\frac{\eta(\tau)\eta(4\tau)}{\eta(2\tau)}.
\]
Their theta series are, respectively,
\[
\sum_{u\in\mathbb Z}(-1)^uq^{u^2},\quad
\frac12\sum_{u\text{ odd}}q^{u^2/8},\quad
\sum_{u\in\mathbb Z}q^{u^2},\quad
\frac12\sum_{u\text{ odd}}\left(\frac8u\right)q^{u^2/8},
\]
where both sums over odd $u$ range over $\mathbb Z$. These give the signs $1,1,1,-1$ at $p=3$. Substituting the three coefficient identities into \eqref{equation_local_hecke_coefficients} with $\lambda=0$ proves \eqref{equation_eta_hecke_eigenvalue} for every primitive quotient.

Finally, suppose $g=V_vg_0$, where $g_0$ is primitive and $V_vF(\tau)=F(v\tau)$. Then $p\nmid v$ and $\Delta_g=v\Delta_{g_0}$, since the weight is $1/2$.
Writing $T_{\nu_g}(p^2)$ and $T_{\nu_{g_0}}(p^2)$ to specify the multipliers, we have
\[
T_{\nu_g}(p^2)V_v=V_vT_{\nu_{g_0}}(p^2).
\]
Indeed, rescaling multiplies the Fourier index by $v$, so the middle Legendre symbol acquires the square factor $v^2$. The first and last terms commute with rescaling because $p\nmid v$. Applying this identity to $g_0$ proves \eqref{equation_eta_hecke_eigenvalue} for $g$ with the same $\chi(p)$.
\end{proof}

\begin{lemma}\label{lemma_eta_local_character}
Under the hypotheses of Lemma~\ref{lemma_eta_hecke_eigenform}, at every cusp with $R_c>0$ we have
\[
\chi(p)=\left(\frac{(d_p/8)\Delta d_ph_c\rho_c}{p}\right),\qquad p\nmid d_p\rho_c.
\]
For $p\ge5$ this symbol is $(\frac{3\Delta R_c}{p})$; for $p=3$ it is $(\frac{\Delta(R_c/3)}3)$. In particular, the expansion at infinity gives the following formula for $\chi(p)$:
\[
\chi(p)=
\begin{cases}
1,&R=0,\\
\displaystyle\left(\frac{(d_p/8)\Delta(d_pR/24)}{p}\right),&R>0,
\end{cases}
\]
and $p\nmid d_pR/24$ when $R>0$.
\end{lemma}
\begin{proof}
Choose $u_p$ as in Theorem~\ref{theorem_hecke_operator}. Since $p\nmid N$, replacing $c$ by $p^2c$ or $u_p^2c$ leaves $h_c,R_c$ unchanged. The least exponent at all these cusps is $\rho_c>0$, so the last term of the Hecke formula at index $d_p\rho_c$ is zero. Comparing with the eigenvalue, we therefore obtain
\[
\frac{a_{p^2c}(p^2d_p\rho_c)}{C_c}-\chi(p)=\frac{\chi(p)-(\frac{(d_p/8)\Delta d_ph_c\rho_c}{p})}{p}.
\]
By Lemma~\ref{lemma_cc_magnitude}, $C_{p^2c}/C_c$ is a root of unity, and $a_{p^2c}(p^2d_p\rho_c)/C_{p^2c}$ is an algebraic integer. The left side is therefore an algebraic integer. The right side is rational, so is an integer; its numerator has absolute value at most $2<p$. The asserted equality and nondivisibility follow.

For the representative $1/N$ of the cusp at infinity, we have $h_N=1$, $R_N=R$, and $\rho_N=R/24$. Thus the assertion for $R>0$ is a special case of the local formula. If $R=0$, the constant term of $g$ at infinity is $1$, and the constant-term formula for the weight $1/2$ Hecke operator shows $\chi(p)=1$.
\end{proof}

\begin{proof}[Proof of Theorem~\ref{theorem_harmonic_weak_maass}]
Use \eqref{equation_uniform_local_indices} with $\kappa=3$ and $\lambda=1$. Write $\xi_{3/2}f=c_0g$, where $c_0$ is an arbitrary constant. By Lemma~\ref{lemma_eta_hecke_eigenform}, the form $g$ is a Hecke eigenform with eigenvalue $(1+\ell^{-1})\chi(\ell)$. By Lemma~\ref{lemma_eta_local_character}, this sign $\chi(\ell)$ equals the value specified in the statement. By the Hecke formula and the definition of $\xi$, we have
\[
\xi_{3/2}(T_{\ell^2}f)=\ell T_{\ell^2}(\xi_{3/2}f)=(1+\ell)\chi(\ell)c_0g.
\]
Here the operator on the shadow has weight $1/2$ and multiplier $\nu_g$. By the defining recurrence and induction,
\[
T_{\ell^{2m}}f-\chi(\ell)^m(1+\ell+\cdots+\ell^m)f\in M_{3/2}^!(\Gamma_0(N),\overline\nu_g).
\]
For $m\ge1$, let
\[
F_m=T_{\ell^{2m}}f-\chi(\ell)T_{\ell^{2m-2}}f-\ell^m\chi(\ell)^mf.
\]
Then $\xi_{3/2}(F_m)=0$, so
\[
F_m\in M_{3/2}^!(\Gamma_0(N),\overline\nu_g),\qquad
H_m:=gF_m\in M_2^!(\Gamma_0(N)).
\]
Theorem~\ref{theorem_hecke_operator} and Lemma~\ref{lemma_cc_magnitude} show that
the Fourier coefficients of $H_m$ at every cusp belong to $\mathcal O$.

At every cusp, (\ref{hypothesis_harmonic_principal}) and Lemma~\ref{lemma_twisted_hecke_principal} give
\[
\mathcal P_c^-(F_m)=\ell^m\sum_{n<0}
\left(a_{u_\ell^{2m}c}(n)q_c^{\ell^{2m}n/d}
-\chi(\ell)^ma_c(n)q_c^{n/d}\right)
\equiv0\pmod{\ell^m}.
\]
The exponents in $F_m|_{3/2}\sigma_c$ belong to $\mathbb Z-\rho_c$. Multiplication by $g_c$ therefore shows that the coefficient of $q_c^n$ in $H_m|_2\sigma_c$ belongs to $\ell^m\mathcal O$ for every integer $n<\lceil\rho_c\rceil$.

We next treat the cusps with $R_c=0$. Write $\chi=\chi(\ell)$, and let $C_c$ denote the nonzero constant term of $g_c=g|_{1/2}\sigma_c$ at such a cusp. By the constant-term formula for the weight $1/2$ Hecke operator and the eigenform relation $T(\ell^2)g=\chi(1+\ell^{-1})g$, we find
\[
C_{\ell^2c}+\ell^{-1}C_{u_\ell^2c}=\chi(1+\ell^{-1})C_c.
\]
Since $\ell u_\ell$ is coprime to $N$, Lemma~\ref{lemma_cc_magnitude} implies
\[
|C_{\ell^2c}|=|C_{u_\ell^2c}|=|C_c|.
\]
Thus, setting $U=C_{\ell^2c}/C_c$ and $V=C_{u_\ell^2c}/C_c$, we have
$|U|=|V|=1$ and
\[
\ell U+V=(\ell+1)\chi.
\]
As $\chi\in\{1,-1\}$, equality in the triangle inequality yields
\begin{equation}\label{equation_zero_cusp_g_constants}
C_{\ell^2c}=C_{u_\ell^2c}=\chi C_c.
\end{equation}

To apply hypothesis (\ref{hypothesis_harmonic_stability}) to an arbitrary representative $1/c$, write $c_* =\gcd(c,N)$ and put
\begin{equation}\label{equation_cusp_stability_modulus}
D_c=\frac{\gcd(c^2,N)}{\gcd(c,N)}=\gcd(c_*,N/c_*).
\end{equation}
The orders in \eqref{equation_explicit_eta_orders} satisfy $\rho_c=\rho_{c_*}$. Since $c_*\mid N$, condition \eqref{equation_main_cusp_stability} therefore gives
$\ell^2\equiv1\pmod{D_c}$ whenever $\rho_c=0$.

We next transfer this relation to forms with the conjugate multiplier. Let $v$ be either $\ell$ or $u_\ell$. By \eqref{equation_main_cusp_stability} and $\ell u_\ell\equiv1\pmod{dN}$, we have $D_c\mid v^2-1$ and $\gcd(v,N)=1$. Thus the congruence
\[
jc^2v^2\equiv c(1-v^2)\pmod N
\]
has an integer solution $j$. For any such solution, put
\[
\gamma_{v,j}=
\begin{pmatrix}
1-jcv^2 & j\\
c(1-v^2-jcv^2) & 1+jc
\end{pmatrix}.
\]
Then $\gamma_{v,j}\in\Gamma_0(N)$ and
\[
\gamma_{v,j}\sigma_{v^2c}=\sigma_c T^j,
\qquad T=\begin{pmatrix}1&1\\0&1\end{pmatrix}.
\]
In particular, the two cusps are equivalent. There is a constant $\mu_{v,c}$, incorporating the multiplier and the cocycle factor arising from the slash operator, such that
\[
g_{v^2c}(\tau)=\mu_{v,c}g_c(\tau+j).
\]
Thus $\mu_{v,c}=\chi$ for both $v=\ell$ and $v=u_\ell$.

For a form $h$ of weight $3/2$ with multiplier $\overline\nu_g$, the product $hg$ has weight $2$ and trivial multiplier. Consequently,
\[
h_{v^2c}(\tau)g_{v^2c}(\tau)=h_c(\tau+j)g_c(\tau+j).
\]
Since an eta-quotient does not vanish on the upper half-plane, we get
\[
h_{v^2c}(\tau)=\mu_{v,c}^{-1}h_c(\tau+j).
\]
Translation does not change the constant coefficient of the holomorphic part. If $b_c(h)$ denotes that coefficient, then
\begin{equation}\label{equation_zero_cusp_h_constants}
b_{\ell^2c}(h)=b_{u_\ell^2c}(h)=\chi b_c(h).
\end{equation}
This relation depends only on the transformation law of $h$ and therefore applies to every Hecke iterate of $f$.

By the weight $3/2$ case of \eqref{equation_local_hecke_coefficients} in Theorem~\ref{theorem_hecke_operator}, we obtain
\[
b_c(T_{\ell^2}h)=b_{\ell^2c}(h)+\ell b_{u_\ell^2c}(h)=(1+\ell)\chi b_c(h).
\]
Set $B_m=b_c(T_{\ell^{2m}}f)$, with $T_1=I$. By \eqref{equation_hecke_recurrence},
\[
B_m=(1+\ell)\chi B_{m-1}-\ell B_{m-2},\qquad
B_0=a_c(0),\quad B_1=(1+\ell)\chi a_c(0).
\]
It follows that
\[
B_m=\chi^m(1+\ell+\cdots+\ell^m)a_c(0).
\]
Hence, for $F_m$ as defined above,
\[
b_c(F_m)=B_m-\chi B_{m-1}-\ell^m\chi^m a_c(0)=0.
\]
Since the principal part of $F_m$ is divisible by $\ell^m$ and $g_c\in\mathcal O[[q_c]]$, the constant term of $H_m=gF_m$ is congruent to $C_cb_c(F_m)=0$ modulo $\ell^m$. Its negative coefficients are also divisible by $\ell^m$. Thus the required bound at these cusps is $b_{\mathfrak a}=1$.

Set $b_{\mathfrak a}=\max\{1,\lceil\rho_{\mathfrak a}\rceil\}$. We have proved that the coefficient of $q_{\mathfrak a}^n$ in $H_m|_2\sigma_{\mathfrak a}$ belongs to $\ell^m\mathcal O$ for every $n<b_{\mathfrak a}$. Moreover,
\[
\sum_{\mathfrak a}(b_{\mathfrak a}-1)=\mathcal W_N(g)\ge2d_2(N)-1.
\]
Lemma~\ref{lemma_integral_cusp_vanishing} therefore gives $H_m\equiv0\pmod{\ell^m}$.
By Lemma~\ref{lemma_cc_magnitude}, we may divide by $g_c$ at every cusp to obtain
\[
F_m\equiv0\pmod{\ell^m}.
\]
The series defining $c_m(n)$ in Lemma~\ref{lemma_fourier_coeffs_fm_infty} is $F_m^++\ell^m\chi(\ell)^mf^+$. Thus $c_m(n)\equiv0\pmod{\ell^m}$ for every $m\ge1$ and every integer $n$.
The two congruences now follow from Lemma~\ref{lemma_fourier_coeffs_fm_infty} as in the proof of Theorem~\ref{theorem_weakly_holomorphic}.
\end{proof}

\begin{proof}[Proof of Remark~\ref{remark_harmonic_constant_alternative}]
Suppose \eqref{equation_harmonic_constant_alternative} holds. The periodicity $f_{c+dN}=f_c$ and $\ell^2u_\ell^2\equiv1\pmod{dN}$ imply
\[
a_{u_\ell^2c}(0)=\chi(\ell) a_c(0)\qquad(\rho_c=0).
\]
For a form $h$ whose constant terms satisfy $b_{\ell^2c}(h)=b_{u_\ell^2c}(h)=\chi(\ell) b_c(h)$, the Hecke coefficient formula gives
\[
b_c(T_{\ell^2}h)=b_{\ell^2c}(h)+\chi(\ell) b_{u_\ell^2c}(h)
=(1+\ell(\ell))\chi b_c(h).
\]
This new vector satisfies the same relations. Hence all Hecke iterates of $f$ do so, and the recurrence in the preceding proof again shows
\[
b_c(T_{\ell^{2m}}f)=\chi(\ell)^m(1+\ell+\cdots+\ell^m)a_c(0).
\]
Thus $b_c(F_m)=0$ at every cusp with $\rho_c=0$. The rest of the proof of Theorem~\ref{theorem_harmonic_weak_maass} uses no instance of (\ref{hypothesis_harmonic_stability}) and applies unchanged.
\end{proof}

\begin{proof}[Proof of Corollary~\ref{corollary_weight_two_dimension_zero}]
The inequality \eqref{equation_main_cusp_inequality} follows from $\mathcal W_N(g)\ge0$ and $d_2(N)=0$. For the genus-zero levels $N=1,\ldots,10,12,13,16,18,25$, every divisor $c\mid N$ satisfies $D_c=\gcd(c,N/c)\in\{1,2,3,4,5\}$. If $D_c\le4$, then $\ell\nmid N$ gives $\ell^2\equiv1\pmod{D_c}$. The case $D_c=5$ occurs only for $N=25$ and $c=5$. Since $r_1+r_5+r_{25}=1$, formula \eqref{equation_explicit_eta_orders} gives
\[
24\rho_5=r_1+5r_5+r_{25}=1+4r_5\ne0.
\]
Thus \eqref{equation_main_cusp_stability} holds at every cusp with $\rho_c=0$.
\end{proof}

\begin{proof}[Proof of Corollary~\ref{corollary_eisenstein_over_eta}]
An eta-quotient has no zeros or poles in $\mathbb H$, and its multiplier has absolute value one. Thus $F=E/g$ belongs to $M_{3/2}^!(\Gamma_0(N),\overline\nu_g)$ and satisfies $\xi_{3/2}F=0$. The product expansion of $g$ gives $F\in q^{-R/24}\mathcal O[[q]]$ at infinity, so Lemma~\ref{lemma_cc_magnitude} gives $\ell$-integrality at every cusp. Since $E$ is holomorphic at $1/c$, the expansion of $g$ there gives
\[
(F|_{3/2}\sigma_c)(\tau)=q_c^{-\rho_c}\bigl(b_c(-d\rho_c)+O(q_c)\bigr).
\]
Since $0\leq\rho_c\leq1$, the only possible polar term is $b_c(-d\rho_c)q_c^{-\rho_c}$, and the principal part vanishes when $\rho_c=0$. Lemma~\ref{lemma_eta_local_character} verifies (\ref{hypothesis_harmonic_principal}) for this polar term. The assertion follows from Corollary~\ref{corollary_weight_two_dimension_zero}.
\end{proof}

\begin{remark}\label{remark_exact_hecke_relation}
Let $g$ be a holomorphic eta-quotient of weight $1/2$ on $\Gamma_0(N)$ satisfying \eqref{equation_main_cusp_inequality}, and let $\ell$ satisfy \eqref{equation_prime_hypotheses} and hypothesis (\ref{hypothesis_harmonic_stability}) of Theorem~\ref{theorem_harmonic_weak_maass}. Let $\chi(\ell)\in\{1,-1\}$ be the sign in \eqref{equation_eta_hecke_eigenvalue}. Suppose that
\[
F\in H_{3/2}^{!}(\Gamma_0(N),\overline{\nu_g}),
\qquad \xi_{3/2}F=c_0g.
\]
Then
\[
T(\ell^2)F=\chi(\ell)(\ell+1)F
\]
if and only if the holomorphic part of $F$ has no negative powers at any cusp.

Suppose that the principal parts vanish, and let $D_\ell=T(\ell^2)F-\chi(\ell)(\ell+1)F$. The identity \eqref{equation_hecke_xi_compatibility} and the coefficient formula of Theorem~\ref{theorem_hecke_operator}, together with \eqref{equation_eta_hecke_eigenvalue}, show that $D_\ell$ is holomorphic on the upper half-plane and at every cusp. The product $gD_\ell$ vanishes at every cusp where $g$ has positive order. At a cusp where $g$ has order zero, the constant-term identity \eqref{equation_zero_cusp_h_constants} and the $n=0$ case of \eqref{equation_local_hecke_coefficients} give $b_c(T(\ell^2)F)=(1+\ell)\chi(\ell)b_c(F)$, so the same conclusion holds. More precisely, its order is at least $b_{\mathfrak a}=\max\{1,\lceil\rho_{\mathfrak a}\rceil\}$ at every cusp. If $gD_\ell$ were nonzero, the valence formula, accounting for forced elliptic zeros, would give
\[
\sum_{\mathfrak a}(b_{\mathfrak a}-1)\le2d_2(N)-2,
\]
contrary to \eqref{equation_main_cusp_inequality}. Hence $gD_\ell=0$.

Conversely, if a principal part is nonzero, choose a polar term at $1/c$ with maximal pole order among all cusps. At $1/(\ell^2c)$, the last term of the Hecke coefficient formula produces a nonzero polar term of $\ell^2$ times that order. Neither the other two terms nor $\chi(\ell)(\ell+1)F$ can have such a pole, so the asserted eigenform relation cannot hold.

This exact relation is analogous to the classical Hecke relation for Zagier’s completed Hurwitz class-number series. Theorem~\ref{theorem_harmonic_weak_maass} also permits prescribed nonzero principal parts, in which case the corresponding Hecke differences need not vanish. Its arithmetic assertion controls these differences modulo powers of a prime through their expansions at all cusps.
\end{remark}

\section{Proofs of the applications}\label{section_application_proofs}
We use $L_t(\tau)=tE_2(t\tau)-E_2(\tau)$ and the notation of
Section~\ref{section_applications}.
\subsection{Sums of core sizes}\label{subsection_proof_core_sizes}

\begin{proof}[Proof of Corollary~\ref{corollary_core_size_congruences}]
Write $F_t(\tau)=\sum_{n\ge0}A_t(n)q^{n-1/24}$, and take $g=\eta$ and $E=tL_t$. Both $E$ and $F_t$ have integral coefficients at infinity, so $F_t$ has $\ell$-integral coefficients at every cusp by Lemma~\ref{lemma_cc_magnitude}. At the cusp $1/c$, put
$v=\gcd(c,t)$ and choose integers $x,y$ with $(t/v)y-(c/v)x=1$. The transformation formula for $E_2$ gives
\begin{equation}\label{equation_lt_cusp_expansion}
(L_t|_2\sigma_c)(\tau)=\frac{v^2}{t}E_2\left(\frac{v^2}{t}\tau-\frac{xv}{t}\right)
-E_2(\tau).
\end{equation}

If $t\le18$, then $d_2(t)=0$ and every cusp width is less than $24$. Thus all cusp orders of $\eta$ are less than $1$, and Corollary~\ref{corollary_eisenstein_over_eta} applies to $A_t(n)$.

For the remaining values of $t$, the cusp $0$ has width $t$. With $S=\begin{pmatrix}0&-1\\1&0\end{pmatrix}$ and $q_0=e^{2\pi i\tau/t}$, we have
\[
(F_t|_{3/2}S)(\tau)=e(1/8)\frac{E_2(\tau/t)-tE_2(\tau)}{\eta(\tau)}.
\]
Its principal part is
\[
e(1/8)q_0^{-t/24}\left(1-t-24\sum_{1\le j<t/24}\sigma_1(j)q_0^j\right).
\]
We call the positive integer $D$ in a term $q_0^{-D/24}$ its polar index. The corresponding Fourier index is $-D$. Here $D=t-24j$ with $0\le j<t/24$, so the polar indices are
\[
\begin{array}{c|rrrrrr}
t&25&27&32&36&49&50\\\hline
D=t-24j&25,1&27,3&32,8&36,12&49,25,1&50,26,2.
\end{array}
\]
At this cusp the width is $t$. For $g=\eta$ we have $\Delta=1$ and
$\chi(\ell)=\left(\frac3\ell\right)$, so the condition
$\epsilon_\ell(h_cD)=\chi(\ell)$ in (\ref{hypothesis_harmonic_principal}) becomes
\[
\left(\frac{tD}{\ell}\right)=1
\quad\text{for each polar index }D.
\]
For $t=25,27,32$, this condition is automatic. For $t=36,50$, it reduces to $\left(\frac3\ell\right)=1$ and $\left(\frac{13}{\ell}\right)=1$, respectively. For $t=49$, the products $tD$ are squares, and $D=25$ requires $\ell\ne5$.

At every other cusp of width $h\le24$, the only possible polar index is $D=h$. The conditions $\ell\nmid D$ and $\epsilon_\ell(hD)=\chi(\ell)$ are automatic, since $\ell\nmid h$ and $\left(\frac{3h^2}{\ell}\right)=\left(\frac3\ell\right)$. The only remaining case is the cusp $1/2$ when $t=50$, whose width is $25$. Here $v=2$, and \eqref{equation_lt_cusp_expansion} shows that the coefficient of $q_c$ in the numerator is zero. The reciprocal eta factor has the form $Cq_c^{-25/24}(1+O(q_c^{25}))$ with $C\ne0$, so the possible polar term $q_c^{-1/24}$ is absent. Thus only the polar index $D=25$ occurs, and $\epsilon_\ell(h_cD)=\left(\frac{3\cdot25^2}{\ell}\right)
=\left(\frac3\ell\right)$. These conditions hold for every equivalent cusp representative, since the polar coefficients change only by roots of unity.
For $t=25$ the global inequality is automatic because $d_2(25)=0$. For $t=27,32,36,49,50$, the divisor formula in Remark~\ref{remark_computable_cusp_criterion} gives
\[
(\mathcal W_{27}(\eta),\mathcal W_{32}(\eta),\mathcal W_{36}(\eta),
\mathcal W_{49}(\eta),\mathcal W_{50}(\eta))=(1,1,1,2,3).
\]
The corresponding values of $d_2(t)$ are $1,1,1,1,2$, so \eqref{equation_main_cusp_inequality} holds in each case. The shadow of $F_t$ is zero, and (\ref{hypothesis_harmonic_stability}) is vacuous since $\eta$ has positive order everywhere. Theorem~\ref{theorem_harmonic_weak_maass} now gives both congruences for $A_t(n)$.

We now pass to $D_t(n)$. By \eqref{equation_prime_level_eisenstein_family},
\[
A_t(n)=D_t(n)+(t-1)(24n-1)p(n).
\]
For the vanishing congruence, put $M=\ell^{2m}n-(\ell^{2m}-1)/24$. Then
$24M-1=\ell^{2m}(24n-1)$, so the correction term is divisible by $\ell^m$. For the recurrence, the two indices have the form $M_r=(\ell^ru+1)/24$, with $r=j$ and $r=j+2$. Since $24M_r-1=\ell^ru$ and $j\ge\lceil j/2\rceil$, the correction
terms are divisible by $\ell^{\lceil j/2\rceil}$ at both indices. Thus the correction can be removed in both congruences.
\end{proof}

\begin{proof}[Proof of Corollary~\ref{corollary_uncorrected_core_congruences}]
The result follows from Corollary~\ref{corollary_core_size_congruences} and
the classical prime-power partition congruences of Watson and Atkin
\cite{watson1938ramanujans,atkin1967proof}, which allow us to remove the partition term at the required indices. 
\end{proof}

\subsection{A second moment of Frobenius partitions}\label{subsection_proof_frobenius}
\begin{proof}[Proof of Corollary~\ref{corollary_frobenius_second_moment}]
Write $\Phi(q)=\prod_{j\ge1}(1-q^j)$ and $c\phi_9(n)=\sum_t c\phi_9(n,t)$. Let $z$ mark the color difference. To select pairs of rows of equal length, we take the constant term in the variable that records the difference of their lengths. Jacobi's triple product then gives
\begin{equation}\label{equation_frobenius_refined_gf}
\sum_{n\ge0}\sum_{t\in\mathbb Z}c\phi_9(n,t)z^tq^n=\frac1{\Phi(q)^9}
\sum_{\substack{\mathbf d\in\mathbb Z^9\\\sum_jd_j=0}}
z^{\sum_{j=0}^8jd_j}q^{\frac12\sum_{j=0}^8d_j^2}.
\end{equation}
At $z=1$ the numerator is the theta series $\Theta_9$ of the root lattice $A_8$. This series belongs to $M_4(\Gamma_0(9))$ \cite{chan2019modular}.

Let $D=q\,d/dq$. To extract the second moment, apply $(z\,\partial/\partial z)^2$ to
\eqref{equation_frobenius_refined_gf} and set $z=1$. By permutation symmetry and the relation $\sum_jd_j=0$, the resulting theta numerator is $15D\Theta_9$. The two generating functions are therefore
\begin{equation}\label{equation_frobenius_theta_gfs}
\sum_{n\ge0}c\phi_9(n)q^n=\frac{\Theta_9}{\Phi^9},\qquad
\sum_{n\ge0}\mathcal M_{2,9}(n)q^n=15\frac{D\Theta_9}{\Phi^9}.
\end{equation}
The second identity also shows that $\mathcal M_{2,9}(n)$ is divisible by $15$, since $A_8$ is an even lattice.

For a modular form $H$ of weight $k$, write
\[
\partial_k H=DH-\frac{k}{12}E_2H
\]
for its Serre derivative. Since $\Theta_9\in M_4(\Gamma_0(9))$, we have $\partial_4\Theta_9\in M_6(\Gamma_0(9))$. We differentiate the first identity in \eqref{equation_frobenius_theta_gfs} and use $D\log\Phi=(E_2-1)/24$ to write the result in the form
\begin{equation}\label{equation_frobenius_corrected_gf}
\frac{27\partial_4\Theta_9}{\eta^9}=\sum_{n\ge0}\left((24n-9)c\phi_9(n)+\frac{\mathcal M_{2,9}(n)}5\right)q^{n-3/8}.
\end{equation}

Take $g=\eta^9$ and $F=27\partial_4\Theta_9/g$. The standard dimension formula gives $d_6(9)=3$ and $d_2(9)=0$. The cusps $\infty,0,1/3,1/6$ have widths $1,9,1,1$, respectively, and
\[
\begin{array}{c|cccc}
\mathfrak a&\infty&0&1/3&1/6\\\hline
\operatorname{ord}_{\mathfrak a}g&3/8&27/8&3/8&3/8.
\end{array}
\]
Hence $\mathcal W_9(g)=3=d_6(9)+d_2(9)$. At the three cusps where $g$ has order $3/8$, holomorphy of $\partial_4\Theta_9$ allows at most a single polar term of exponent $-3/8$. It remains to examine the cusp $0$. By \cite[Theorem~5.1]{chan2019modular},
\begin{equation}\label{equation_frobenius_eta_expansion}
\begin{aligned}
\Theta_9={}&324\eta(3\tau)^8
+19683\frac{\eta(9\tau)^{12}}{\eta(3\tau)^4}
-240\frac{\eta(\tau)^9\eta(9\tau)^3}{\eta(3\tau)^4}\\
&\quad-1458\frac{\eta(\tau)^6\eta(9\tau)^6}{\eta(3\tau)^4}
+\frac{\eta(\tau)^{12}}{\eta(3\tau)^4}.
\end{aligned}
\end{equation}
We use the transformation law for $\eta$ on each of the five terms. In the local parameter $t=e(\tau/9)$, the expansion at $0$ is
\[
\Theta_9|_4S=\frac13+6t^4+O(t^7).
\]
Since the Serre derivative commutes with the slash action, the expansion of $F$ at $0$ is
\[
(F|_{3/2}S)(\tau)=e(1/8)t^{-27/8}\bigl(-3+O(t^4)\bigr).
\]
Since $4>27/8$, only the leading term contributes to the principal part.

The numerator $H=27\partial_4\Theta_9$ has constant term $-9$ at infinity, whereas $H|_6S=-3+O(t^4)$. If $H$ were Eisenstein, its vanishing coefficients of $t,t^2,t^3$ at $0$ would force it to be a multiple of $E_6(\tau)$, as one sees from the standard Eisenstein basis of $M_6(\Gamma_0(9))$. Such a multiple has equal constant terms at $0$ and infinity. Hence $H$ has a nonzero cusp-form component.

The coefficients of $F$ at infinity are integers by \eqref{equation_frobenius_corrected_gf}, since $15\mid\mathcal M_{2,9}(n)$.

Here $r=R=9$, $\Delta=1$, $\kappa=3$, $k=6$, and $\lambda=1$. At every cusp, $h_c(24\rho_c)=9h_c^2$, so $\epsilon_\ell(dh_c\rho_c)=\left(\frac3\ell\right)=\chi(\ell)$, and the polar index is prime to $\ell\ge5$. Theorem~\ref{theorem_weakly_holomorphic} applies to
$A_9(n)=(24n-9)c\phi_9(n)+\mathcal M_{2,9}(n)/5$. The condition in part (\ref{conclusion_weak_zero}) is $\left(\frac{9-24n}{\ell}\right)=1$. For the vanishing congruence, put $M=\ell^{2m}n-3(\ell^{2m}-1)/8$. Then $24M-9=\ell^{2m}(24n-9)$, so the correction term $(24M-9)c\phi_9(M)$ is divisible by $\ell^m$. For the recurrence, take $n=3u$ in Theorem~\ref{theorem_weakly_holomorphic}. The corresponding indices are $M_r=(\ell^ru+3)/8$, with $r=j$ and $r=j+2$, and $24M_r-9=3\ell^ru$. Thus $(24M_r-9)c\phi_9(M_r)$ is divisible by $\ell^{\lceil j/2\rceil}$ at both indices. The two congruences for $\mathcal M_{2,9}(n)/5$ hold modulo $\ell^m$ and $\ell^{\lceil j/2\rceil}$, respectively. For the moment itself, each modulus acquires a factor of $5$. Since $\ell\ge5$, these moduli are prime to $3$, and \eqref{equation_frobenius_theta_gfs} shows that $3$ divides every value of the moment. This accounts for the remaining factor of $3$.
\end{proof}

\subsection{Further weakly holomorphic applications}
\begin{proof}[Proof of Corollary~\ref{corollary_overpartition_cubes}]
Let
\[
g(\tau)=\frac{\eta(\tau)^2}{\eta(2\tau)},\qquad
E(\tau)=\frac{E_4(\tau)-E_4(2\tau)}{120},\qquad F=E/g.
\]
Then $E\in M_4(\Gamma_0(2))$, and $d_4(2)=d_2(2)=0$, so \eqref{equation_dimension_cusp_criterion} is automatic. The orders of $g$ at $\infty$ and $0$ are $0$ and $1/8$, respectively. By \eqref{equation_overpartition_cubes_gf}, $F$ has integral coefficients and zero constant term at infinity. The constant term of $E|_4S$ is $(1-2^{-4})/120=1/128$, so $F$ has a single polar term of exponent $-1/8$ at $0$.

Here $R=0$, $\Delta=1/2$, $\kappa=7$, and $\lambda=3$. At the only polar cusp, $h_0=2$ and $\rho_0=1/8$, so
\[
\epsilon_\ell(d_\ell h_0\rho_0)=\left(\frac{(d_\ell/8)^2}{\ell}\right)=1.
\]
Theorem~\ref{theorem_weakly_holomorphic} therefore applies with $\chi(\ell)=1$. The indices are unshifted because $R=0$, and $\epsilon_\ell(-d_\ell n)=\left(\frac{-n}{\ell}\right)$. The two congruences have moduli $\ell^{5m}$ and $\ell^{5\lceil j/2\rceil}$, respectively.
\end{proof}

\begin{proof}[Proof of Corollary~\ref{corollary_odd_three_color}]
To prove \eqref{equation_odd_three_color_gf}, weight each odd part of size $j$ by $z^j$ in the product for $P_3$ and differentiate at $z=1$. Since
\[
L_2(\tau)=1+24\sum_{\substack{j\ge1\\j\text{ odd}}}\frac{jq^j}{1-q^j},\qquad
\frac1{g_3(\tau)}=q^{1/24}P_3(q),
\]
we obtain \eqref{equation_odd_three_color_quotient}. The cusp orders of $g_3$ on $\Gamma_0(2)$ are $\rho_\infty=-1/24$ and $\rho_0=1/6$, with widths $1$ and $2$,
respectively. By \eqref{equation_odd_three_color_quotient}, $L_2/g_3$ has integral coefficients at infinity and vanishes there. At $0$, the transformation formulas give
\[
g_3|_{1/2}S=2e(-1/8)\frac{\eta(\tau)^3}{\eta(\tau/2)^2},
\qquad L_2|_2S=\tfrac12E_2(\tau/2)-E_2(\tau).
\]
In the parameter $q_0=e(\tau/2)$, the principal part of $L_2/g_3$ is $-e(1/8)q_0^{-1/6}/4$. With $r=1$, $k=2$, $R=-1$, and $\Delta=1/4$, we have
\[
\mathcal W_2(g_3)=-1=2d_2(2)-1,\qquad \epsilon_\ell(d_\ell h_0\rho_0)=\left(\frac6\ell\right).
\]
The two congruences for $b_3$ now follow from Theorem~\ref{theorem_weakly_holomorphic} for every $\ell\ge5$. In part (\ref{conclusion_weak_zero}) the sign condition reads
$\left(\frac{-2(24n+1)}\ell\right)=1$.

For $\ell=7$, Chen \cite[Theorem~1.5]{chen2025congruences} proved that $c_3(n)\equiv0\pmod{7^{\lceil a/2\rceil}}$ whenever $a\ge1$ and $7^a\mid24n+1$. At the indices of the two congruences, the $c_3$ terms therefore vanish modulo $7^m$ and $7^{\lceil j/2\rceil}$, respectively. Since $8$ is prime to $7$, the identity $b_3=c_3+8O_3$ transfers both congruences to $O_3$. The condition $\left(\frac{-2(24n+1)}7\right)=1$ is equivalent to $n\equiv3,4,6\pmod7$.
\end{proof}

\subsection{Applications of Theorem \ref{theorem_harmonic_weak_maass}}\label{subsection_proof_harmonic_applications}
\begin{proof}[Proof of Corollary~\ref{corollary_sodo_congruences}] Identity \eqref{equation_sodo_eisenstein} follows by marking the odd parts in the generating function $Q(q)$.

We have $d_2(4)=0$, and the cusp widths and orders are
\[
\begin{array}{c|ccc}
c&1&2&4\\\hline
h_c&4&1&1\\
\operatorname{ord}_{1/c}(g_O)&1/8&0&1/8.
\end{array}
\]
The numerator $2L_2-L_4$ is a weight $2$ modular form on $\Gamma_0(4)$ with integral coefficients at infinity. By \eqref{equation_lt_cusp_expansion}, its constant term at $1/2$ is $2$, whereas $g_O$ has order zero. Corollary~\ref{corollary_eisenstein_over_eta}, the zero-shadow case of the harmonic theorem, now applies at every odd prime with $R=3$, $\Delta=2$, and $\chi(\ell)=\left(\frac2\ell\right)$.
\end{proof}

\subsubsection{Smallest-parts functions}\label{subsection_proof_smallest_parts}
Let $\mathcal A$ and $\mathcal G$ be the completed forms in \cite[Sections 1 and 2]{ahlgren2015euler}, normalized by
\[
\begin{aligned}
\mathcal A^+(\tau)&
=\sum_{n\ge1}\spt(n)q^{n-1/24}-\frac{E_2(\tau)}{12\eta(\tau)},\\
\mathcal G^+(\tau)
&=\frac{2\eta(2\tau)}{\eta(\tau)^2}\left(\sum_{n\ge1}\frac{(-1)^nq^{n(n+1)}}{(1-q^n)^2}-\frac1{24}-\frac{E_2(2\tau)}{12}\right).
\end{aligned}
\]
Put
\[
\begin{aligned}
g_1(\tau)&=\eta(2\tau),&
f_1(\tau)&=\mathcal A(2\tau)+\frac{L_2(\tau)}{24g_1(\tau)},\\
g_2(\tau)&=\frac{\eta(2\tau)^2}{\eta(4\tau)},&
f_2(\tau)&=\mathcal G(2\tau)+\frac{L_4(\tau)}{24g_2(\tau)}.
\end{aligned}
\]
By the generating-function identities used in \cite[Corollaries 19 and 20]{wang2026explicit}, the holomorphic parts have the expansions
\[
\begin{aligned}
f_1^+(\tau)&=\sum_{n\ge0}
\left(\spt_\omega(n)+\frac{12n-1}{24}p(n/2)\right)q^{n-1/12},\\
f_2^+(\tau)&=\sum_{n\ge0}
\left(\overline{\spt}_\omega(n)-\frac18\overline p(n/2)\right)q^n.
\end{aligned}
\]

\begin{lemma}\label{lemma_application_cusp_data}
Use $\gamma_\infty=I$, $\gamma_0=S=\left(\begin{smallmatrix}0&-1\\1&0\end{smallmatrix}\right)$, and $\gamma_{1/2}=U=\left(\begin{smallmatrix}1&0\\2&1\end{smallmatrix}\right)$. For $t=e(\tau/h_{\mathfrak a})$, write
$g_i|_{1/2}\gamma_{\mathfrak a}=Ct^\rho(1+O(t))$. Let $P$ be the sum of the terms of exponent at most zero in $(f_i|_{3/2}\gamma_{\mathfrak a})^+$, and let $b^-(0)$ be the coefficient of $y^{-1/2}$ in its non-holomorphic part. Then
\begingroup
\renewcommand{\arraystretch}{1.65}
\[
\begin{array}{c|c|c|c|c|c|c}
i&\mathfrak a&h_{\mathfrak a}&\rho&C&P&b^-(0)\\ \hline
1&\infty&1&1/12&1&-\frac1{24}t^{-1/12}&0\\
1&0&2&1/24&\frac{e(-1/8)}{\sqrt2}
&-\frac{e(1/8)}{12\sqrt2}t^{-1/24}&0\\ \hline
2&\infty&1&0&1&-\frac18&\frac1{2\sqrt2\pi}\\
2&0&4&1/8&e(-1/8)&-\frac{e(1/8)}{32}t^{-1/8}&0\\
2&1/2&1&1/8&\sqrt2 e(-1/16)&0&0
\end{array}
\]
\endgroup
All coefficients of $(f_i|_{3/2}\gamma_{\mathfrak a})^+$ belong to
$\mathbb Z[\zeta_{48},1/6]$, and for $i=2$ they belong to $\mathbb Z[\zeta_{48},1/2]$. Moreover,
\[
\xi_{3/2}f_1=-\frac{\sqrt3}{4\pi}g_1,\qquad
\xi_{3/2}f_2=-\frac1{4\sqrt2\pi}g_2.
\]
\end{lemma}

\begin{proof}
We recall from \cite[Sections 1 and 2]{ahlgren2015euler} that
\[
\mathcal A(-1/\tau)=e(1/8)\tau^{3/2}\mathcal A(\tau),\qquad
\mathcal G(-1/(2\tau))=-2(-i\tau)^{3/2}\mathcal H(\tau),
\]
where
\[
\mathcal H^+(\tau)=\sum_{n\ge1}(-1)^n\operatorname{M2spt}(n)q^{n-1/8}+\frac{\eta(\tau)}{24\eta(2\tau)^2}\bigl(E_2(2\tau)-E_2(\tau)\bigr)
\in q^{7/8}\mathbb Z[[q]].
\]
Here $\operatorname{M2spt}$ is the integer-valued statistic in that paper. These identities and the transformation laws for $\eta$ and $E_2$ determine the cusp expansions in the table.
For $i=2$, consider the terms at $0$ and $1/2$ whose denominators could be divisible by $3$. Their respective expansions, with $u=\tau/4$ in the first, are
\[
\begin{aligned}
\frac{\frac14E_2(u)-E_2(4u)}{24}
&=-\frac1{32}-\frac14\sum_{n\ge1}\sigma_1(n)e(nu)+\sum_{n\ge1}\sigma_1(n)e(4nu),\\
\frac{E_2(\tau-1/2)-E_2(\tau)}{24}&=\sum_{n\ge1}(1-(-1)^n)\sigma_1(n)q^n.
\end{aligned}
\]
Neither expansion has a denominator divisible by $3$. Together with $\mathcal H^+\in q^{7/8}\mathbb Z[[q]]$ and the expansion at infinity, this proves the integrality assertions. The shadow formulas follow from the completions in \cite{ahlgren2015euler}. At each cusp the constant term of the shadow determines the last column of the table, since $\xi_{3/2}(by^{-1/2})=-\overline b/2$.
\end{proof}

\begin{proof}[Proof of Corollary~\ref{corollary_spt_omega}]
Lemma~\ref{lemma_application_cusp_data} establishes integrality at the cusps of $\Gamma_0(2)$ and $\Gamma_0(4)$ and identifies the shadows required in (\ref{hypothesis_harmonic_shadow}). The coefficients at any equivalent representative $1/c$ differ only by factors that are roots of unity. Each nonzero principal part consists of a single term of exponent $-\rho_c$, so (\ref{hypothesis_harmonic_principal}) follows from Lemma~\ref{lemma_eta_local_character}. Since $d_2(2)=d_2(4)=0$, Corollary~\ref{corollary_weight_two_dimension_zero} gives \eqref{equation_main_cusp_inequality} and (\ref{hypothesis_harmonic_stability}).

We now apply Theorem~\ref{theorem_harmonic_weak_maass} to $f_1$ with $R=\Delta=2$ and $\ell\ge5$. At the indices in both congruences, the factor $12n-1$ in the partition correction is divisible by the required power of $\ell$, as in the proof of Corollary~\ref{corollary_core_size_congruences}. Since $24$ is a unit at $\ell$, the correction vanishes modulo that power, and the congruences hold for $\spt_\omega$ itself.

For $f_2$, the parameters are $R=0$ and $\Delta=1$, and $f_2^+$ already has the coefficients in the second row of the corollary. At every cusp the holomorphic coefficients lie in $\mathbb Z[\zeta_{48},1/2]$, so the theorem applies at $\ell=3$ as well. In each case the vanishing condition of the theorem reduces to the quadratic-symbol condition in the corollary.
\end{proof}

\section*{AI declaration}
After the earlier version of the manuscript was finished, the author used GPT Astra 6.0 (OpenAI) to review and polish the paper. The original version considered only congruences for $p\geq5$. The case $p=3$, including the Hecke operators and the corresponding congruences, was treated with the help of Astra. The author also asked Astra to find  more examples of the results, and the author carefully chose a few of them. Some of the mathematical exposition, proof simplification, and the strengthening of key conclusions were carried out with the help of Astra. The author independently reviewed and verified all AI-assisted content and takes full responsibility for the results, proofs, and examples presented in this work.

\providecommand{\bysame}{\leavevmode\hbox to3em{\hrulefill}\thinspace}
\providecommand{\MR}{\relax\ifhmode\unskip\space\fi MR }
\providecommand{\MRhref}[2]{%
  \href{http://www.ams.org/mathscinet-getitem?mr=#1}{#2}
}
\providecommand{\href}[2]{#2}


\begin{thebibliography}{10}

\bibitem{ahlgren2024shimura}
S. Ahlgren, N. Andersen, and R. Dicks, \emph{The Shimura lift and
congruences for modular forms with the eta multiplier}, Adv. Math.
\textbf{447} (2024), 109655.

\bibitem{ahlgren2015euler}
S. Ahlgren and N. Andersen, \emph{Euler-like recurrences for smallest parts functions}, Ramanujan J. \textbf{36} (2015), 237--248.


\bibitem{zbMATH05941220}
S. Ahlgren, K. Bringmann, and J. Lovejoy, \emph{{$\ell$}-adic
  properties of smallest parts functions}, Adv. Math. \textbf{228} (2011),
  no.~1, 629--645.

\bibitem{zbMATH06342625}
S. Ahlgren and B. Kim, \emph{Mock modular grids and {Hecke} relations
  for mock modular forms}, Forum Math. \textbf{26} (2014), no.~4, 1261--1287.

\bibitem{ahlgren2009vanishing}
S. Ahlgren, N. Masri, and J. Rouse,
\emph{Vanishing of modular forms at infinity},
Proc. Amer. Math. Soc. \textbf{137} (2009), no.~4, 1205--1214.

\bibitem{amdeberhan2025crank}
T. Amdeberhan and M. Merca, \emph{From crank to congruences},
Mediterr. J. Math. \textbf{22} (2025): 190.

\bibitem{andrews1984generalized}
G.~E. Andrews, \emph{Generalized Frobenius partitions},
Mem. Amer. Math. Soc. \textbf{49} (1984), no.~301.

\bibitem{zbMATH06828531}
G.~E. Andrews, A. Dixit, D. Schultz, and A.~J. Yee,
  \emph{Overpartitions related to the mock theta function {$\omega(q)$}},
  Acta Arith. \textbf{181} (2017), no.~3, 253--286.


\bibitem{zbMATH06533997}
G.~E. Andrews, A. Dixit, and A.~J. Yee, \emph{Partitions associated with
  the {Ramanujan}/{Watson} mock theta functions {$\omega(q)$},
  {$\nu(q)$} and {$\phi(q)$}}, Res. Number Theory \textbf{1} (2015): 19.

\bibitem{atkin1967proof}
A.~O.~L. Atkin, \emph{Proof of a conjecture of Ramanujan}, Glasgow Math. J.
\textbf{8} (1967), 14--32.



\bibitem{bardhan2026proofs}
G. Bardhan and N. Saikia, \emph{Proofs of the conjectures on
{$SOME(n)$} and {$DSOME(n)$} functions related to integer partitions},
preprint (2026), arXiv:2608.08524 [math.NT].

\bibitem{bhattacharya2017holomorphic}
S. Bhattacharya, \emph{Holomorphic eta quotients of weight $1/2$}, Adv. Math. \textbf{320} (2017), 1185--1200.


\bibitem{zbMATH06828732}
K. Bringmann, A. Folsom, K. Ono, and L. Rolen, \emph{Harmonic
  {Maass} forms and mock modular forms: theory and applications}, American Mathematical Society Colloquium Publications, vol.~64, American Mathematical Society, Providence, RI, 2017.
  
    



\bibitem{bruinier2004two}
J.~H. Bruinier and J. Funke, \emph{On two geometric theta lifts}, Duke
  Math. J. \textbf{125} (2004), no.~1, 45--90.

\bibitem{chan2019modular}
H.~H. Chan, L. Wang, and Y. Yang, \emph{Modular forms and $k$-colored
generalized Frobenius partitions}, Trans. Amer. Math. Soc.
\textbf{371} (2019), no.~3, 2159--2205.

\bibitem{chen2025congruences}
D. Chen, \emph{Congruences modulo powers of $7$ for the reciprocal
crank parity function}, preprint (2025), arXiv:2511.18905 [math.CO].

\bibitem{cohen2017modular}
H. Cohen and F. Str{\"o}mberg, \emph{Modular forms: a classical
  approach}, Graduate Studies in Mathematics, vol. 179, American Mathematical Society, Providence, RI, 2017.


\bibitem{diamond2005first}
F. Diamond and J. Shurman, \emph{A first course in modular forms},
Graduate Texts in Mathematics, vol.~228, Springer, New York, 2005.

\bibitem{garvan2026sum}
F. Garvan and R. Sarma, \emph{Sum of parts in overpartitions and
partitions without repeated odd parts}, preprint (2026),
arXiv:2606.14680 [math.NT].

\bibitem{katz1973p}
N. M. Katz, \emph{$p$-adic properties of modular schemes and modular forms},
in Modular Functions of One Variable III, Lecture Notes in Math.
\textbf{350}, Springer, Berlin, 1973, pp. 69--190.

\bibitem{knopp2008modular}
M.~I. Knopp, \emph{Modular functions in analytic number theory}, AMS Chelsea Publishing, vol. 337, American Mathematical Society, Providence, RI, 2008.

\bibitem{kolitsch1992generalized}
L.~W. Kolitsch, \emph{Generalized Frobenius partitions, $k$-cores,
$k$-quotients, and cranks}, Acta Arith. \textbf{62} (1992), no.~1,
97--102.

\bibitem{lemkeoliver2013eta}
R.~J. Lemke Oliver, \emph{Eta-quotients and theta functions}, Adv.
  Math. \textbf{241} (2013), 1--17.

\bibitem{zbMATH07567828}
R. Mao, \emph{Congruences for
  {$\overline{\operatorname{spt}}_\omega(n)$} modulo powers of primes}, J.
  Number Theory \textbf{240} (2022), 107--123.

\bibitem{zbMATH08135742}
\bysame, \emph{Congruences for the smallest parts function associated with
  {$\omega (q)$}}, J. Comb. Theory Ser. A \textbf{217} (2026), 106106.

\bibitem{zbMATH00963239}
Y. Martin, \emph{Multiplicative {$\eta$}-quotients}, Trans. Amer. Math.
  Soc. \textbf{348} (1996), no.~12, 4825--4856.



\bibitem{merca2024sums}
M. Merca, \emph{From sums of divisors to partition congruences}, Rev. R. Acad.
Cienc. Exactas F\'{i}s. Nat. Ser. A Mat. RACSAM \textbf{118} (2024),
Article 115.


\bibitem{zbMATH06133929}
K. Ono, \emph{Congruences for the {Andrews} spt function}, Proc. Natl. Acad.
  Sci. USA \textbf{108} (2011), no.~2, 473--476.



\bibitem{rolen2017t}
L. Rolen, \emph{On $t$-core towers and $t$-defects of partitions},
Ann. Comb. \textbf{21} (2017), no.~1, 119--130.

\bibitem{sturm1987congruence}
J. Sturm, \emph{On the congruence of modular forms}, in Number Theory (New York, 1984--1985), Lecture Notes in Math. \textbf{1240}, Springer, Berlin, 1987, 275--280.


  
\bibitem{zbMATH07595066}
L. Wang and Y. Yang, \emph{Moments of ranks and cranks, and quotients
  of {Eisenstein} series and the {Dedekind} eta function}, J. Number Theory
  \textbf{242} (2023), 471--521.
  
\bibitem{wang2026explicit}
W. Wang, \emph{Explicit constructions for {Ramanujan}-type congruences}, preprint (2026), arXiv:2606.19008 [math.NT]. 



\bibitem{watson1938ramanujans}
G.~N. Watson, \emph{Ramanujans Vermutung \"uber Zerf\"allungsanzahlen},
J. Reine Angew. Math. \textbf{179} (1938), 97--128.

\bibitem{williams2025hecke}
C. Williams, \emph{Hecke relations for eta multipliers and congruences for
  higher-order smallest parts functions}, J. Number Theory \textbf{277}
  (2025), 325--343.




\bibitem{zhou2023double}
H.-G. Zhou and X.-J. Zhu, \emph{Double coset operators and eta-quotients},
J. Number Theory \textbf{249} (2023), 537--601.

\bibitem{zhu2024dimension}
X.-J. Zhu, \emph{Dimension formulas for modular form spaces of rational
weights, the classification of eta-quotient characters and an extension
of Martin's theorem}, preprint (2024), arXiv:2408.00246 [math.NT].

\end{thebibliography}
\end{document}